\def\date{24 March 2017} 
\magnification \magstep1
\overfullrule=0pt
\newif\ifproofmode
\def\xrefsfilename{c5c.xrf}  

\newcount\referno
\newcount\thmno
\newcount\secno
\referno=0
\thmno=0
\secno=0
\def\ifundefined#1{\expandafter\ifx\csname#1\endcsname\relax}
\expandafter \def \csname SECLABELintro\endcsname {1}
\expandafter \def \csname THMLABEL1.1\endcsname {(1.1)}
\expandafter \def \csname THMLABEL1.2\endcsname {(1.2)}
\expandafter \def \csname THMLABEL1.3\endcsname {(1.3)}
\expandafter \def \csname THMLABEL1.4\endcsname {(1.4)}
\expandafter \def \csname THMLABEL1.5\endcsname {(1.5)}
\expandafter \def \csname THMLABEL1.6\endcsname {(1.6)}
\expandafter \def \csname THMLABEL1.7\endcsname {(1.7)}
\expandafter \def \csname SECLABELext\endcsname {2}
\expandafter \def \csname THMLABEL1extquadconn\endcsname {(2.1)}
\expandafter \def \csname THMLABELlongext\endcsname {(2.2)}
\expandafter \def \csname THMLABEL2.2\endcsname {(2.3)}
\expandafter \def \csname THMLABEL2.3\endcsname {(2.4)}
\expandafter \def \csname THMLABEL2.4\endcsname {(2.5)}
\expandafter \def \csname THMLABEL2.6\endcsname {(2.6)}
\expandafter \def \csname SECLABELhe\endcsname {3}
\expandafter \def \csname THMLABEL3.3\endcsname {(3.1)}
\expandafter \def \csname THMLABEL3.4\endcsname {(3.2)}
\expandafter \def \csname THMLABEL3.5\endcsname {(3.3)}
\expandafter \def \csname THMLABEL3.6\endcsname {(3.4)}
\expandafter \def \csname SECLABELquad\endcsname {4}
\expandafter \def \csname THMLABELextex\endcsname {(4.1)}
\expandafter \def \csname THMLABELezremark\endcsname {(4.2)}
\expandafter \def \csname THMLABEL4.1\endcsname {(4.3)}
\expandafter \def \csname THMLABEL4.2\endcsname {(4.4)}
\expandafter \def \csname THMLABEL4.3\endcsname {(4.5)}
\expandafter \def \csname THMLABEL4.4\endcsname {(4.6)}
\expandafter \def \csname THMLABEL4.5\endcsname {(4.7)}
\expandafter \def \csname THMLABEL4.55\endcsname {(4.8)}
\expandafter \def \csname THMLABEL4.6\endcsname {(4.9)}
\expandafter \def \csname THMLABEL4.6cor\endcsname {(4.10)}
\expandafter \def \csname SECLABELdodec\endcsname {5}
\expandafter \def \csname THMLABEL5.1\endcsname {(5.1)}
\expandafter \def \csname THMLABEL5.2\endcsname {(5.2)}
\expandafter \def \csname THMLABEL5.3\endcsname {(5.3)}
\expandafter \def \csname THMLABEL5.35\endcsname {(5.4)}
\expandafter \def \csname THMLABEL5.4\endcsname {(5.5)}
\expandafter \def \csname THMLABEL5.5\endcsname {(5.6)}
\expandafter \def \csname SECLABELtwo\endcsname {6}
\expandafter \def \csname THMLABEL6.1\endcsname {(6.1)}
\expandafter \def \csname THMLABEL6.2\endcsname {(6.2)}
\expandafter \def \csname SECLABELone\endcsname {7}
\expandafter \def \csname THMLABEL7.1\endcsname {(7.1)}
\expandafter \def \csname THMLABEL7.2\endcsname {(7.2)}
\expandafter \def \csname THMLABEL7.3\endcsname {(7.3)}
\expandafter \def \csname THMLABEL7.4\endcsname {(7.4)}
\expandafter \def \csname THMLABEL7.5\endcsname {(7.5)}
\expandafter \def \csname THMLABEL7.6\endcsname {(7.6)}
\expandafter \def \csname REFLABELAldHolJac\endcsname {1}
\expandafter \def \csname REFLABELBar\endcsname {2}
\expandafter \def \csname REFLABELBut\endcsname {3}
\expandafter \def \csname REFLABELDieGT\endcsname {4}
\expandafter \def \csname REFLABELEdwSanSeyTho\endcsname {5}
\expandafter \def \csname REFLABELMcCPhD\endcsname {6}
\expandafter \def \csname REFLABELMcCEdge\endcsname {7}
\expandafter \def \csname REFLABELRobSanSeyTho4CT\endcsname {8}
\expandafter \def \csname REFLABELRobSeyGM9\endcsname {9}
\expandafter \def \csname REFLABELRobSeyThoTut\endcsname {10}
\expandafter \def \csname REFLABELRobSeyThoCubic\endcsname {11}
\expandafter \def \csname REFLABELTutConvex\endcsname {12}
\expandafter \def \csname REFLABELTutGeom\endcsname {13}

\immediate\openout1=\xrefsfilename
\def\bibitem#1#2\par{\ifundefined{REFLABEL#1}\relax\else
 \global\advance\referno by 1\relax
 \immediate\write1{\noexpand\expandafter\noexpand\def
 \noexpand\csname REFLABEL#1\endcsname{\the\referno}}
 \global\expandafter\edef\csname REFLABEL#1\endcsname{\the\referno}
 \item{\the\referno.}#2\ifproofmode [#1]\fi\fi}
\def\cite#1{\ifundefined{REFLABEL#1}\ignorespaces
   \global\expandafter\edef\csname REFLABEL#1\endcsname{?}\ignorespaces
   \write16{ ***Undefined reference #1*** }\fi
 \csname REFLABEL#1\endcsname}
\def\nocite#1{\ifundefined{REFLABEL#1}\ignorespaces
   \global\expandafter\edef\csname REFLABEL#1\endcsname{?}\ignorespaces
   \write16{ ***Undefined reference #1*** }\fi}
\def\newthm#1#2\par{\global\advance\thmno by 1\relax
 \immediate\write1{\noexpand\expandafter\noexpand\def
 \noexpand\csname THMLABEL#1\endcsname{(\the\secno.\the\thmno)}}
 \global\expandafter\edef\csname THMLABEL#1\endcsname{(\the\secno.\the\thmno)}
 \medbreak\noindent{\bf(\the\secno.\the\thmno)\enspace}\ignorespaces
 \ifproofmode {\bf[#1]} \fi{\sl#2}
 \ifdim\lastskip<\medskipamount\removelastskip\penalty55\medskip\fi}
\def\newsection#1#2\par{\global\advance\secno by 1\relax
 \immediate\write1{\noexpand\expandafter\noexpand\def
 \noexpand\csname SECLABEL#1\endcsname{\the\secno}}
 \global\expandafter\edef\csname SECLABEL#1\endcsname{\the\secno}
 \vskip0pt plus.3\vsize
 \vskip0pt plus-.3\vsize\bigskip\bigskip\vskip\parskip
 \message{\the\secno. #2}\thmno=0
 \centerline{\bf\the\secno. #2\ifproofmode {\rm[#1]} \fi}
 \nobreak\smallskip\noindent}
\def\refthm#1{\ifundefined{THMLABEL#1}\ignorespaces
 \global\expandafter\edef\csname THMLABEL#1\endcsname{(?)}\ignorespaces
 \write16{ ***Undefined theorem label #1*** }\fi
 {\csname THMLABEL#1\endcsname}}
\def\refsec#1{\ifundefined{SECLABEL#1}\ignorespaces
 \global\expandafter\edef\csname SECLABEL#1\endcsname{(?)}\ignorespaces
 \write16{ ***Undefined section label #1*** }\fi
 \csname SECLABEL#1\endcsname}
\outer\def\thm#1#2\par{\medbreak\noindent{\bf(#1)\enspace}\ignorespaces
{\sl#2}\ifdim\lastskip<\medskipamount\removelastskip\penalty55\medskip\fi}
\def\cond#1#2\par{\smallbreak\noindent\rlap{\rm(#1)}\ignorespaces
\hangindent=36pt\hskip36pt{\rm#2}\smallskip}
\def\claim#1#2\par{\smallbreak\noindent\rlap{\rm(#1)}\ignorespaces
\hangindent=30pt\hskip30pt{\sl#2}\smallskip}
\def\dfn#1{{\sl #1}}

\def\proof{\smallbreak\noindent{\sl Proof. }}



\def\qed{\hfill$\square$\medskip}
\def\sqr#1#2{{\vcenter{\vbox{\hrule height.#2pt
\hbox{\vrule width.#2pt height #1pt \kern#1pt
\vrule width.#2pt}
\hrule height.#2pt}}}}
\def\square{\mathchoice\sqr56\sqr56\sqr{2.1}3\sqr{1.5}3}
\def\restriction{|}
\def\ref#1#2{\item{#1.}#2}
\outer\def\beginsection#1\par{\vskip0pt plus.3\vsize
   \vskip0pt plus-.3\vsize\bigskip\bigskip\vskip\parskip
   \message{#1}\centerline{\bf#1}\nobreak\smallskip\noindent}
\input colordvi

\footline={\ifnum\pageno=1\hfil\else\rm\hfil\folio\hfil\fi}
\def\junk#1{}
\def\emb{\hookrightarrow}
\def\cfc{cyclically $5$-con\-nected}
\def\cfcity{cyclic $5$-con\-nectivity}
\def\dc{dodecahedrally connected}
\def\dcity{dodecahedral con\-nectivity}
\def\he{homeomorphic embedding}
\def\qc{quad-con\-nected}
\def\qcity{quad-con\-nectivity}
\hyphenation{trip-lex}
\font\smallrm=cmr8

\phantom{a}
\bigskip\bigskip
\centerline{\bf CYCLICALLY FIVE--CONNECTED CUBIC GRAPHS}
\bigskip\bigskip\bigskip
\baselineskip=11pt

\centerline{Neil Robertson$^{*1}$\vfootnote{$^*$}{\smallrm
Research partially 
performed under a consulting agreement with Bellcore, and partially
supported by DIMACS Center, 
Rutgers University, New Brunswick, New Jersey  08903, USA.}
\vfootnote{$^1$}{\smallrm Partially supported
by NSF under Grant No. DMS-8903132 and by ONR under Grant No. 
N00014-92-J-1965.
}}
\centerline{Department of Mathematics}
\centerline{Ohio State University}
\centerline{231 W. 18th Ave.}
\centerline{Columbus, Ohio  43210, USA}
\bigskip

\centerline{P. D. Seymour$^{2}$\vfootnote{$^2$}{\smallrm
This research was performed while the author was employed at Bellcore, 445 South St., Morristown, NJ 07960.
}}
\centerline{Department of Mathematics}
\centerline{Princeton University}
\centerline{Princeton, New Jersey  08544, USA}
\bigskip

\centerline{and}
\bigskip

\centerline{Robin Thomas$^{*3}$\vfootnote{$^3$}{\smallrm
Partially supported
by NSF under Grants No. DMS-9303761 and  DMS-1202640
and by ONR under Grant No. 
N00014-93-1-0325.
}}
\centerline{School of Mathematics}
\centerline{Georgia Institute of Technology}
\centerline{Atlanta, Georgia  30332, USA}
\vfill
\noindent 23 February 1995

\noindent Revised  \date
\vfil
\eject

\footline{\hss\tenrm\folio\hss}
\baselineskip 18pt
\beginsection ABSTRACT

\bigskip

\noindent
A cubic graph $G$ is \dfn{\cfc} if $G$ is simple, $3$--connected,
has at least $10$ vertices
and for every set $F$ of edges of 
size at most four, at most one component of $G\backslash F$ contains
circuits. We prove that if $G$ and $H$ are \cfc\ cubic graphs
and $H$ topologically contains $G$, then either $G$ and $H$ are
isomorphic, or (modulo well-described exceptions) there exists
a \cfc\ cubic graph $G'$ such that $H$ topologically contains $G'$
and $G'$ is obtained from $G$ in one of the following  two ways.
Either $G'$ is obtained from $G$ by subdividing two
distinct edges of $G$ and joining the two new vertices by an edge,
or $G'$ is obtained from $G$ by subdividing each edge of a circuit of
length five and joining the new vertices by a matching to a new circuit
of length five disjoint from $G$ in such a way that the cyclic orders
of the two circuits agree. We prove a companion result, where
by slightly increasing the connectivity of $H$ we are able to eliminate
the second construction.
We also prove versions of both of these results when $G$ is almost
\cfc\ in the sense that it satisfies the definition except for
$4$-edge cuts such that one side is a circuit of length four.
In this case $G'$ is required to be almost \cfc\ and to have
fewer circuits of length four than $G$. In particular,
if $G$ has at most one circuit of length four, then
$G'$ is required to be \cfc. However, in this more general
setting the
operations  describing the possible graphs $G' $  are more complicated. 

\vfil\eject

\newsection{intro} INTRODUCTION

The primary motivation for this work comes from Tutte's 
$3$--edge-coloring conjecture [\cite{TutGeom}], the following 
(definitions are given later).

\newthm{1.1}{\bf Conjecture.} Every $2$--edge-connected cubic graph that does not 
topologically contain the Petersen graph is $3$--edge-colorable.

\noindent
Our strategy is to reduce (1.1) to ``apex" and ``doublecross" graphs,
two classes of graphs that are close to planar graphs, and then modify
our proof of the Four Color Theorem~[\cite{RobSanSeyTho4CT}] to
show that graphs belonging to those classes satisfy (1.1).
We began the first part of this program in~[\cite{RobSeyThoTut}],
but in order to complete it we need to understand the structure
of reasonably well-connected cubic graphs that do not topologically
contain the Petersen graph.
That is the subject of~[\cite{RobSeyThoCubic}], where we
apply the structure theory of \cfc\ cubic graphs developed in this paper.
We have completed the second part of the project for doublecross
graphs in~[\cite{EdwSanSeyTho}]; 
the apex case is harder and is currently under preparation.

To motivate our structure theorems let us mention a special case
of a theorem of Tutte~[\cite{TutConvex}]. 

\newthm{1.2}Let $G,H$ be non-isomorphic $3$--connected cubic graphs, 
and let $H$ contain $G$ topologically. 
Then there exists a cubic graph $G'$ obtained from $G$ 
by subdividing two distinct edges of $G$ and joining the
new vertices by an edge in such a way such that $H$ 
topologically contains $G'$.

Our objective is to prove a similar theorem for \cfc\ cubic graphs.
An ideal analog of (1.2) for \cfc\ cubic graphs would assert that
there is a graph $G'$ as in (1.2) that is \cfc. 
That is unfortunately not true, but the exceptions can
be conveniently described. We will do so now.

Let $G$ be a \cfc\ cubic graph. Let $e,f$ be distinct edges of $G$
with no common end and such that no edge of $G$ is adjacent
to both $e$ and $f$,
and let $G'$ be obtained from $G$ by subdividing $e$ and $f$  
and joining the new vertices by an edge. We say that $G'$ is 
a \dfn{handle expansion} of $G$. 
We show in \refthm{longext} that $G'$ is \cfc.
Let $e_1,e_2,e_3,e_4,e_5$ 
(in order) be the 
edges of a circuit of $G$ of length five. Let us subdivide
$e_i$ by a new vertex $v_i$, add a circuit (disjoint from $G$)
with vertices $u_1,u_2,u_3,u_4,u_5$ (in order), and for 
$i=1,2,\ldots,5$ let us add an edge joining $u_i$ and $v_i$ to
form a graph $G''$. In these circumstances we say
that $G''$ is a \dfn{circuit expansion} of $G$.
It is not hard to see, for instance by repeatedly
applying~\refthm{1extquadconn},  that $G''$ is \cfc.

Let $p$ be an integer such that $p\ge5$ if $p$ is odd and $p\ge10$
if $p$ is even. Let $G$ be a cubic graph with vertex-set 
$\{u_0,u_1,\ldots,u_{p-1},v_0,v_1,\ldots,v_{p-1}\}$ such that
for $i=0,1,\ldots,p-1$, 
$u_i$ has neighbors $u_{i-1}$, $u_{i+1}$ and $v_i$, and $v_i$
has neighbors $u_i$, $v_{i-2}$ and $v_{i+2}$, where the index
arithmetic is taken modulo $p$ (see Figure 1). We say that $G$ is
a \dfn{biladder} on $2p$ vertices. We remark that the Petersen
graph is a biladder on $10$ vertices, and that the Dodecahedron
is a biladder on $20$ vertices.

\goodbreak\midinsert
\vskip4.2in
\includegraphics{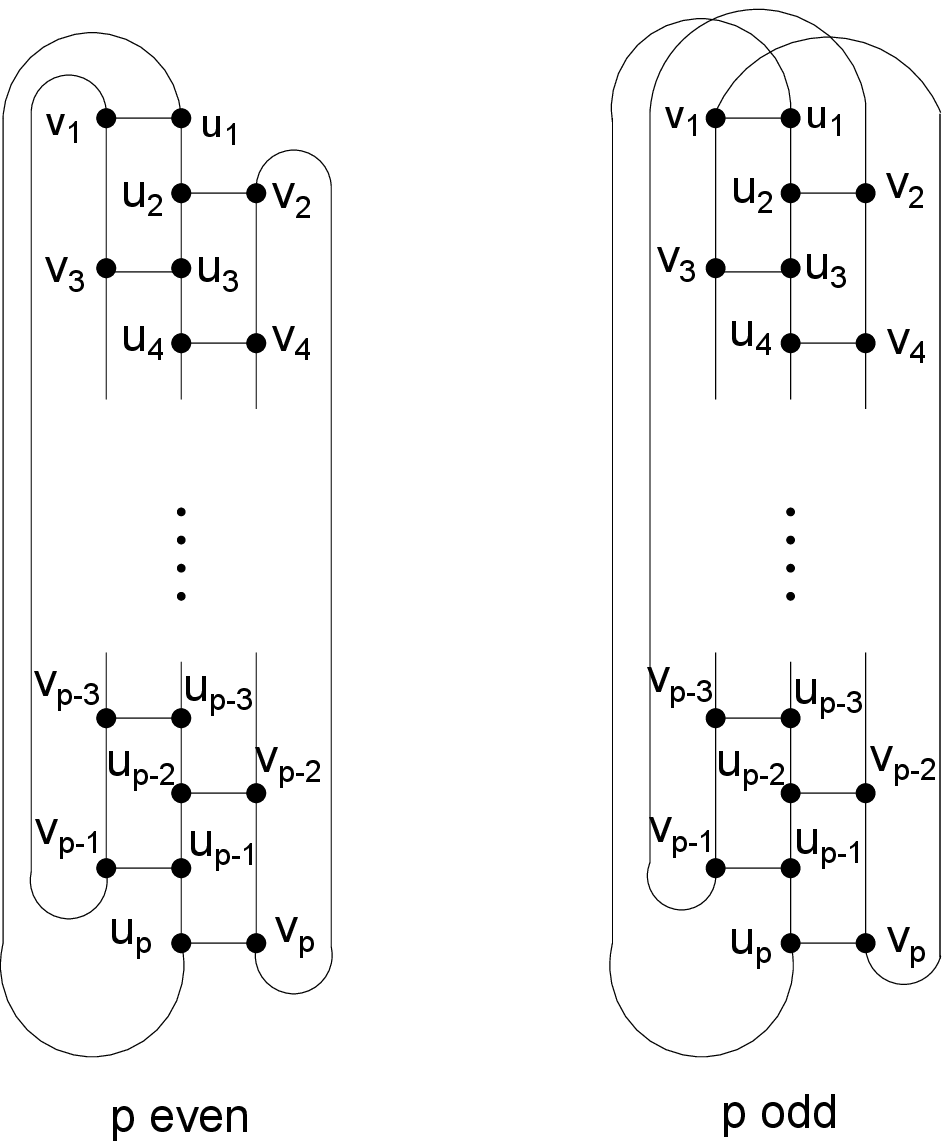}
\centerline{Figure 1: Biladders}
\endinsert


\junk{Let $n\ge5$ be an integer. Let $C$ be a circuit on the vertices
$v_1,v_2,\dots,v_{2n}$ (in order), let $C_1$ be a circuit on 
$u_1,u_3,\ldots,u_{2n-1}$, and let $C_2$ be a circuit on
$u_2,u_4,\ldots,u_{2n}$.
The \dfn{biladder}
on $4n$ vertices is the cubic graph $G$ with vertex-set
$\{u_1,u_2,\ldots,u_{2n},v_1,v_2,\ldots,v_{2n}\}$ obtained from the 
disjoint union of $C,C_1$ and $C_2$ by adding an edge between
$u_i$ and $v_i$ for all $i=1,2,\ldots,2n$. We remark that the
biladder on $20$ vertices is isomorphic to the Dodecahedron, and that
every biladder on at least $20$ vertices is \cfc.
Let $m\ge2$ be an integer. Let $C'$ be a circuit on the vertices
$v_1,v_2,\dots,\allowbreak v_{2m+1}$ (in order), and
let $C''$ be a circuit on 
$u_1,u_3,\ldots,u_{2m+1},u_2,u_4,\ldots,\allowbreak u_{2m}$.
The \dfn{M\"obius biladder}
on $4m+2$ vertices is the cubic graph $G$ with vertex-set
$\{u_1,u_2,\ldots,u_{2m+1},v_1,v_2,\ldots,v_{2m+1}\}$ 
obtained from the 
disjoint union of $C$ and $C'$ by adding an edge between
$u_i$ and $v_i$ for all $i=1,2,\ldots,2m+1$. We remark that
the M\"obius biladder on $10$ vertices is isomorphic to the Petersen
graph, and that every M\"obius biladder on at least $10$ vertices
is \cfc.}

The following is our first main result.

\newthm{1.3}Let $G,H$ be non-isomorphic \cfc\ cubic graphs such that
not both of them are biladders,  let 
$H$ contain $G$
topologically, and assume
that if $G$ is isomorphic
to the Petersen graph, then $H$ does not topologically contain the
biladder on $14$ vertices (that is, for $p=7$), and
 if $G$ is isomorphic to the Dodecahedron, then $H$ does
not topologically contain the biladder on $24$ vertices (that is, for $p=12$).
Then there exists a \cfc\ handle or circuit expansion $G'$ 
of $G$ such that $H$ contains $G'$ topologically.

There is a variation of (1.3), which is easier to apply, but which
involves a stronger assumption about the graph $H$. Dodecahedral
connectivity is defined in Section 5.

\newthm{1.4}Let $G,H$ be non-isomorphic \cfc\ cubic graphs such that
not both of them are biladders, 
let $H$ be \dc,  let $H$ contain $G$ topologically, and assume
that if $G$ is isomorphic
to the Petersen graph, then $H$ does not topologically contain the
biladder on $14$ vertices (that is, for $p=7$), and
if $G$ is isomorphic to the Dodecahedron, then $H$ does
not topologically contain the biladder on $24$ vertices (that is, for $p=12$).
Then there exists a \cfc\ handle expansion $G'$ of $G$
such that $H$ topologically contains $G'$.

Since every biladder is either planar (if $p$ is even), or
topologically contains the Petersen graph (if $p$ is odd) we
deduce the following corollary.

\newthm{1.5}Let $G,H$ be non-isomorphic \cfc\ cubic graphs,
let $G$ be non-planar,
let $H$ be \dc,  let $H$ contain $G$ topologically, and assume
that $H$ does not topologically contain the Petersen graph.
Then there exists a \cfc\ handle expansion $G'$ of $G$
such that $H$ topologically contains $G'$.

The last three theorems describe how to obtain a bigger \cfc\ cubic
graph from a smaller one. But what are the initial graphs to start
from? The graphs
Petersen, Triplex, Box, Ruby and Dodecahedron are defined in Figure 2.
The following theorem of McCuaig [\cite{McCPhD}, \cite{McCEdge}] was
also obtained in [\cite{AldHolJac}].

\goodbreak\midinsert
\vskip4.2in
\includegraphics{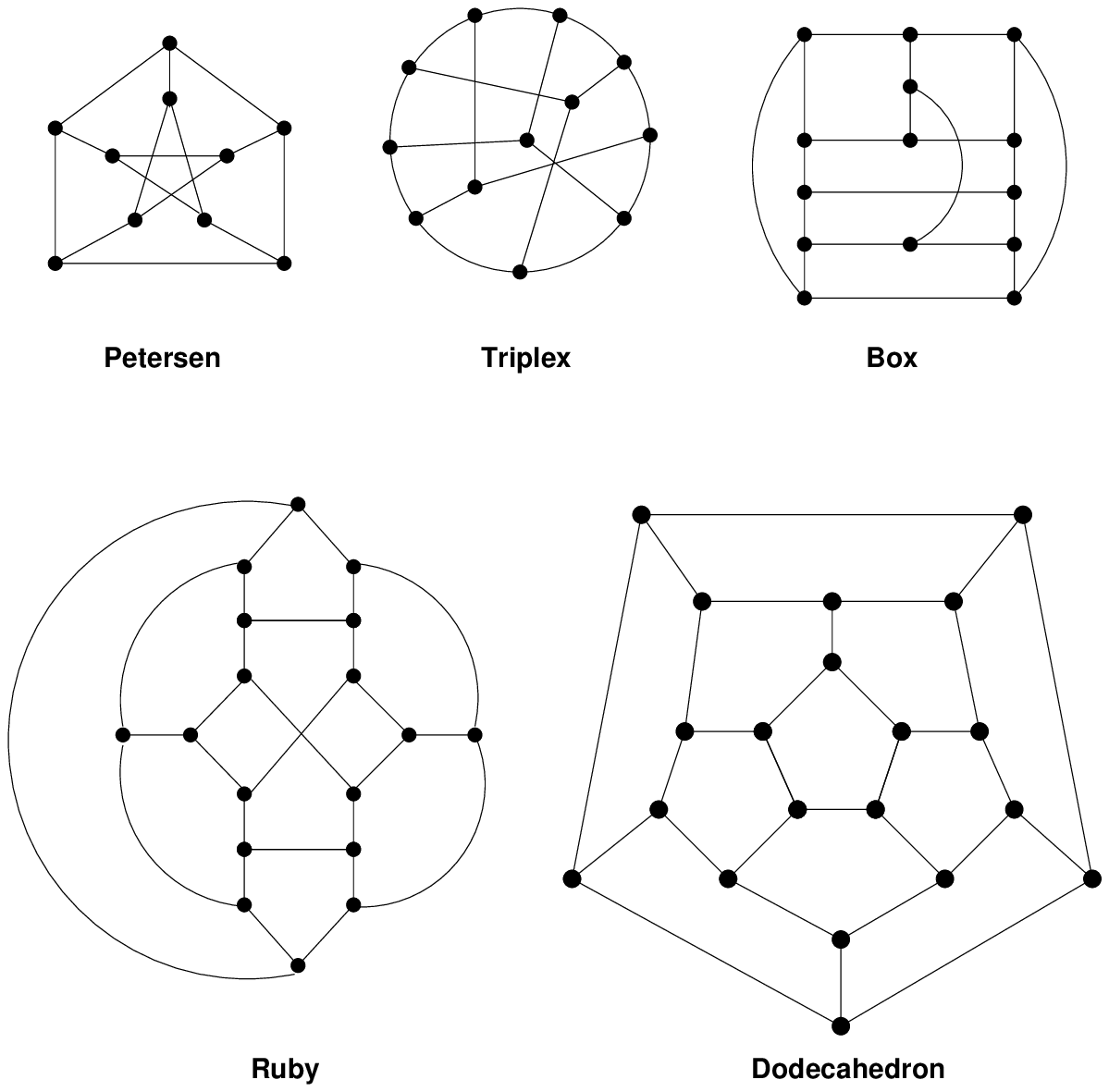}
\centerline{Figure 2. The five minimal \cfc\ graphs}
\endinsert

\newthm{1.6}Every \cfc\ cubic graph topologically contains one of
Petersen, Triplex, Box, Ruby or Dodecahedron.

Theorems (1.3), (1.4) and (1.6) have the following corollary,
the first part of which was
proved for planar graphs in [\cite{Bar}, \cite{But}], and for general
graphs in [\cite{McCPhD}, \cite{McCEdge}].

\newthm{1.7}Every \cfc\ cubic graph can be obtained from
Triplex, Box, Ruby or a biladder 
by repeatedly applying the operations of handle expansion
or circuit expansion. Every \dc\ cubic graph can be
obtained from Triplex, Box, Ruby or a biladder 
by repeatedly applying the operation of handle expansion.

\noindent
It follows from~\refthm{5.1} that a handle expansion of a \dc\
graph is again \dc. 

Our proofs of (1.3) and (1.4) are indirect, and proceed by way
of auxiliary results, stated as \refthm{4.6cor} and \refthm{5.5} below,
that are themselves quite useful.
Those auxiliary results allow $G$ to violate the definition of
\cfcity, but only in a limited way. 
For instance, if $G$ satisfies the definition of \cfcity, except for
one circuit of length four, then we can still insist that $G'$ be \cfc.
However, the operations that describe how to obtain $G'$ are
more complicated, and therefore we defer the exact statements
to Section~\refsec{quad}.

Let us introduce some terminology now. All graphs in this paper
are finite and simple. Thus we may denote
the edge of a graph with ends $u$ and $v$
by $uv$ without any ambiguity.
If $G$ is a graph we denote its vertex-set and edge-set by
$V(G)$ and $E(G)$, respectively. Let $G$ be a graph.
If $K,L$ are subgraphs of $G$ we denote by $K\cup L$ the graph
with vertex-set $V(K)\cup V(L)$, edge-set $E(K)\cup E(L)$ and
the obvious incidences.
If $A\subseteq V(G)$ we denote by $\delta_G A$ (or $\delta A$ if
the graph can be understood from the context)
the set of edges
of $G$ with one end in $A$ and the other end in $V(G)-A$.
An \dfn{edge-cut} of $G$ is a set of edges of the form $\delta A$,
where $A\subseteq V(G)$ and $\emptyset\not=A\not= V(G)$. 
If $X$ is a vertex, a set of vertices, an edge, or a set of edges, we
denote by $G\backslash X$ the graph obtained from $G$ by deleting $X$.
If $X$ is a set of vertices we denote by $G\restriction X$ the 
graph $G\backslash(V(G)-X)$.
\dfn{Paths} and \dfn{circuits} have no ``repeated" vertices and no 
``repeated" edges. A \dfn{quadrangle} is a circuit of length four.
A graph $G$ is \dfn{cubic} if every vertex of $G$ has degree
three and it is \dfn{subcubic} if every vertex has degree at most three.
Let $k\ge4$ be an integer. We say that a cubic graph $G$ is 
\dfn{cyclically $k$--connected} if $G$ is $3$--connected, 
has at least $2k$ vertices, and 
for every edge-cut $\delta A$ of $G$ of cardinality less than $k$,
one of $G\backslash A$, $G\restriction A$ has no circuits.

Let $e$ be an edge of a graph $G$. A graph $H$ is obtained from
$G$ by \dfn{subdividing} $e$ if $H$ is obtained by deleting $e$,
adding a new vertex $v\not\in V(G)$, and joining $v$ to both ends
of $e$ by new edges. We say that $v$ is the \dfn{new vertex} of $H$.
We say that a graph $H$ \dfn{topologically contains} a graph $G$
if some graph obtained from $G$ by repeatedly subdividing edges
is isomorphic to a subgraph of $H$.

The paper is organized as follows. In Sections 2 and 3 we introduce
some terminology and prove several lemmas. In Section 4
we solve the 
following problem: Suppose that a \cfc\ cubic graph $H$
contains a graph $G$ topologically and is minimal with this property, 
where $G$ is  ``almost" \cfc\ (\qc, as defined in the next section). 
 What can we say about $H$? In Section 5 we strengthen
the conclusion of the result of Section 4 under the 
assumption that $H$ is
\dc. In Section 6 we prove a preliminary version of (1.3), where
we allow adding two handles, rather than one.
We prove (1.3) and (1.4)
in Section 7.

\newsection{ext} EXTENSIONS

Let $G$ be a cubic graph.
We say that $G$ is \dfn{\qc} if 
\item{$\bullet$} $G$ is cyclically $4$-connected, 
\item{$\bullet$} $G$ has at least $10$ vertices, 
\item{$\bullet$} if $G$ has more than one quadrangle, 
then it has at least $12$ vertices, and 
\item{$\bullet$} for every edge-cut
$\delta A$ of $G$ of cardinality exactly four, one of
$G\restriction A$ and $G\backslash A$ is isomorphic to $K_2$ or to a quadrangle.

\noindent
Thus a cubic graph is \cfc\ if and only if it is \qc\ and has no quadrangle.
It follows that in a \qc\ graph no two quadrangles share an edge.

Let $u,v,x,y$ be vertices of a graph $G$ such that $u$ is
adjacent to $v$, $x$ is adjacent to $y$, and $\{u,v\}\ne\{x,y\}$.
 We define
$G+(u,v,x,y)$ to be the graph obtained from $G$ by subdividing the edges
$uv$ and $xy$, where the new vertices are $k$ and $l$, respectively,
and adding an edge joining $k$ and $l$. The vertices $k,l$ (in
this order) will be called the \dfn{new vertices} of $G+(u,v,x,y)$.
We remark that if $u$, $v$, $x$, $y$ are pairwise distinct and 
$G$ has no circuits of length at most three,
then neither does $G+(u,v,x,y)$.
If the vertices $u,v,x,y$ are pairwise distinct, then we say that
 $G+(u,v,x,y)$ is a \dfn{$1$-extension} of $G$.
If, in addition, neither $u$ nor $v$ is adjacent to $x$ or $y$, then we say that
 $G+(u,v,x,y)$ is a \dfn{long $1$-extension of $G$};
otherwise we say that it is \dfn{short}.
Thus if $G $ is \cfc, then long $1$-extension and handle expansion mean the same thing.

Now let $C$ be a quadrangle in $G$.
We say that the 1-extension  $G+(u,v,x,y)$ is \dfn{based at $C$}
if $uv$ is an edge of $C$ and $x,y\not\in V(C)$, and we will apply
the qualifiers long and short as in the previous paragraph.
Let $G$ be \qc, let $C$ be a quadrangle in $G$, let 
$G+(u,v,x,y)$ be a short 1-extension of $G$ based at $C$, 
and let $k,l$ be the new vertices.
Then one of $u,v$ is adjacent to one of $x,y$, and so
we may assume that, 
%
%
say, $u$ is adjacent to $x$. Then $\{u,x,k,l\}$ 
is the vertex-set of a quadrangle $D$ in $G'$.
The next lemma implies that $D$ is the only quadrangle of $G'$ containing the edge $kl$.
We say that $D$ is the \dfn{new quadrangle} of $G'$.

\newthm{1extquadconn} Let $G$ be a \qc\ cubic graph, 
and let $G'$ be a $1$-extension of $G$ 
such that if $G$ has a quadrangle, then $G'$ is a $1$-extension of
$G$ based at some quadrangle of $G$.
Then $G'$ is \qc.
In particular, $G'$ has at most one quadrangle that is not a quadrangle of $G$.

\proof 
\junk{
Let $G'=G+(u,v,x,y)$.  We first prove the following.
\claim{1}There is at most one pair $a,b$ of adjacent vertices of
$G$ with $a\in\{u,v\}$ and $b\in\{x,y\}$.
To prove (1) suppose to the contrary that two such pairs exist,
and let $u,x$ be one such pair. Since $G$ has no circuits of
length three we see that the other pair is $v,y$. Then $u,v,y,x$
are the vertices of a quadrangle $D$ in $G$.  
Since $G$ has a quadrangle we see that $G'$ is a 1-extension 
based at some quadrangle $C$, where $u,v\in V(C)$ and $x,y\not\in V(C)$.
Thus $C\ne D$, and yet $C$ and $D$ have at least one edge in common, 
contrary to the \qcity\ of $G$.  This proves (1).
\medskip
}
The second assertion follows from the first,
for every quadrangle of $G'$ that is not a quadrangle of $G$ uses the edge joining
the new vertices of $G'$ and in a \qc\ graph every edge belongs to at most one quadrangle.
To prove that $G'$ is \qc\ it suffices to verify the last condition
in the definition of \qcity, because the other conditions are clear. 
The graph $G'$ is clearly $3$--connected.
Let $k,l$ be the new vertices of $G'$. 
Let $\delta_{G'} A$ be an edge-cut
of $G'$ of cardinality at most four such that both $G'\restriction A$
and $G'\backslash A$ contain circuits. 
We must show that $|\delta_{G'}A|=4$ and that $G'|A$ or
$G'\backslash A$ is a quadrangle.
We have $4\le|A|\le |V(G')|-4$.
Let $B=A-\{k,l\}$. Then $\emptyset\not=B\not=V(G)$, and so
$\delta_G B$ is an edge-cut of $G$ of cardinality at most four. 
Thus one of
$G\restriction B$ and $G\backslash B$ is a forest or a quadrangle.

Suppose first that $G|B$ is a quadrangle.  Since $4=|\delta_GB|\le|\delta_{G'}A|\le4$,
we see that $\delta_GB=\delta_{G'}A$.
The definition of 1-extension implies that $\{u,v\}\not\subseteq B$ and
$\{x,y\}\not\subseteq B$.  
Thus $G|B=G'|A$, and so $G'|A$ is a quadrangle, as desired. 
This completes the case where $G|B$ is a quadrangle.

By symmetry between $G|B$ and $G\backslash B$ we may therefore assume that
$G|B$ is a forest.  Since $|\delta_GB|\le 4$ we see that $|B|\le 2$, and
since $|A|\ge 4$ we have $|B|=2$, say $B=\{a,b\}$.  
Thus $A=\{a,b,k,l\}$, $|\delta_{G'}A|=4$, and $G'|A$ is a quadrangle, as required.\qed

\newthm{longext} Let $G$ be a \qc\ cubic graph with at most one
quadrangle, and let $G'$ be a long $1$-extension such that if
$G$ has a quadrangle $C$, then $G'$ is a $1$-extension based at $C$.
Then $G'$ is \cfc.

\proof
The graph $G'$ is \qc\ by \refthm{1extquadconn}. 
Since the extension is long, the graph $G'$ has no quadrangle,
and hence is \cfc.~\qed 

\newthm{2.2}Let $G$ be a \qc\ cubic graph, 
let the vertices $u_1,u_2,\allowbreak u_3,\allowbreak u_4,u_5$
(in order) form the vertex-set of a path of $G$, let 
$G'=G+(u_1,u_2,u_4,u_5)$, and assume that 
either $G$ is \cfc, or $G$ has a quadrangle $C$
 with $u_1,u_2\in V(C)$ and $u_4,u_5\not\in V(C)$. 
Then $G'$ is a short extension of $G$ if and only if $u_1$ and
$u_5$ are adjacent.

\proof If $u_1$ and $u_5$ are adjacent, 
then $G'$ is clearly a short extension.
Conversely, if $G'$ is a short extension,  
then one of $u_1,u_2$ is adjacent to one of 
$u_4,u_5$. Since $G$ has no triangles we may assume for
a contradiction that either $u_1$ is adjacent to $u_4$,
or $u_2$ is adjacent to $u_5$. In either case $G$ has a
quadrangle $D\ne C$, and hence $G$ is not \cfc.
Thus $C$ exists, but the existence of $C$ and $D$ contradicts
the \qcity\ of $G$.~\qed

Let $G$ be a cyclically $4$--connected cubic graph, let 
$u_1,u_2,\ldots,u_6$ be the vertices of a path in $G$ in order, 
let $G_1=G+(u_1,u_2,u_3,u_4)$, and let $k_1,l_1$ be the new
vertices of $G_1$. We define $G\&(u_1,u_2,u_3,u_4,u_5,u_6)$
to be the graph $G_2=G_1+(u_3,l_1,u_5,u_6)$. Let $k_2,l_2$ be
the new vertices of $G_2$. We say that $k_1,l_1,k_2,l_2$ are the
\dfn{new vertices} of $G\&(u_1,u_2,u_3,u_4,u_5,u_6)$.

\newthm{2.3}Let $G$ be a  cubic graph, and let 
$u_1,u_2,\ldots,u_6$ be vertices of $G$ forming the vertex-set
of a path in the order listed. Let $G'=G\&(u_1,u_2,u_3,u_4,u_5,u_6)$.
Assume that $G$ is \qc\ with at most one quadrangle, 
and that if it has a quadrangle, then it has a quadrangle $C$
with $u_1,u_2\in V(C)$ and $u_4,u_5,u_6\not\in
V(C)$. Then $G'$ is \cfc.

\proof By \refthm{1extquadconn} $G_1=G+(u_1,u_2,u_3,u_4)$ is \qc,
and it has exactly one quadrangle.
By another
application of \refthm{1extquadconn} the graph $G'$ is \cfc, because it has no 
quadrangle.~\qed

\newthm{2.4}Let $G$ be a cubic graph, let $u_1,u_2,\ldots,u_5$ be
the vertices of a circuit of $G$ in order, 
and assume that $G$ is either \cfc, or \qc\ with a  quadrangle
$C$ such that  $u_3,u_4\in V(C)$ and
$u_1,u_2,u_5\not\in V(C)$.
Let $v_1$ be the
neighbor of $u_1$ other than $u_2$ and $u_5$, 
and let $G'=G+(u_3,u_4,u_1,v_1)$.
 Then $G'$ is a long $1$-extension of $G$.

\proof The vertex
$u_3$ is not adjacent to $v_1$ in $G$, for otherwise $G$ has a quadrangle
$D$ with vertex-set $\{v_1,u_1,u_2,u_3\}$, which implies that $C$ exists,
but the existence of $C$ and $D$ contradicts the \qcity\ of $G$. 
Hence $G'$ is a long $1$-extension by \refthm{2.2}. \qed

\newthm{2.6}Let $G$ be a \qc\ cubic graph, and
let $C$ be a quadrangle in $G$ with vertices $u_1,u_2,u_3,u_4$
in order. 
Let $v_1$ be the neighbor of $u_1$
not on $C$, and let $v_2$ be defined analogously. Let $w_1\not=u_1$
be a neighbor of $v_1$, and let $z_1\not=v_1$ be a neighbor of 
$w_1$. Then $G+(u_2,u_3,v_1,w_1)$ is a long $1$-extension of $G$, 
and if $z_1\not=v_2$ then $G+(u_1,u_2,w_1,z_1)$ is a long $1$-extension of $G$.

\proof The vertices $w_1$ and $u_3$ are not adjacent, for otherwise
the set $\{u_1,u_2,u_3,u_4,v_1,w_1\}$ contradicts the \qcity\ of $G$.
Thus the $1$-extension $G+(u_3,u_4,v_1,w_1)$ is long by~\refthm{2.2}, 
and so is 
$G+(u_1,u_2,w_1,z_1)$ if $z_1\not=v_2$.~\qed

\newsection{he} HOMEOMORPHIC EMBEDDINGS

Let $G,H$ be graphs. A mapping $\eta$ with domain $V(G)\cup E(G)$
is called a \dfn{homeomorphic
embedding} of $G$ into $H$ if for every two vertices $v,v'$  and 
every two
edges $e,e'$ of $G$
\item{(i)}$\eta(v)$ is a vertex of $H$, and if $v,v'$ are distinct
then $\eta(v),\eta(v')$ are distinct,
\item{(ii)}if $e$ has ends $v,v'$, then $\eta(e)$ is a path of $H$
with ends $\eta(v),\eta(v')$, and otherwise disjoint from $\eta(V(G))$, and
\item{(iii)}if $e,e'$ are distinct, then $\eta(e)$ and $\eta(e')$ are
edge-disjoint, and if they have a vertex in common, then this vertex
is an end of both.


\noindent We shall denote the fact that $\eta$ is a 
homeomorphic embedding of $G$ into $H$ by writing $\eta:G\emb H$,
and we shall write $G\emb H$ to mean that there exists a \he\
of $G$ into $H$.
If $K$ is a subgraph of $G$ 
we denote by $\eta(K)$ the subgraph of $H$ consisting of all
vertices $\eta(v)$, where $v\in V(K)$, and all vertices and 
edges that belong
to $\eta(e)$ for some $e\in E(K)$.
It is easy to see that $H$ contains $G$ topologically if and only
if there is a homeomorphic embedding $G\emb H$. 

Let $G_0$ be a \qc\ graph,  let $C_0$ be a quadrangle in $G$, and let $n\ge1$
be an integer.
We say that $G_n$ is an \dfn{$n$-extension} of $G_0$ \dfn{based at} $C_0$ if
there exists a sequence $G_1,C_1,G_2, C_2,\dots, G_n$ such that for
$i=1,2,\dots, n$, the graph $G_i$ is a 1-extension of $G_{i-1}$ based at
$C_{i-1}$, and if $i<n$ then this 1-extension is short and $C_i$ is the
new quadrangle in $G_i$.  We say that $G_n$ is a \dfn{short}
\dfn{$n$-extension} of $G_0$ if $G_n$ is a short $1$-extension of $G_{n-1}$, 
and we say that it is a \dfn{long $n$-extension} otherwise.
We say that the sequence $G_1,G_2, \dots, G_n$
 is a \dfn{generating sequence of the $n$-extension $G_n$
from $G_0$ based at $C_0$}.
We say that a graph $H$ is an \dfn{extension} of $G_0$ if it is
an $n$-extension for some integer $n\ge1$. 

Let $G,H,K$ be graphs, and let $\eta:G\emb H$ and $\zeta:H\emb K$.
For $v\in V(G)$ we put $\xi(v)=\zeta(\eta(v))$, and for $e\in
E(G)$ we define $\xi(e)$ to be the union of $\zeta(f)$ over all
edges $f\in E(\eta(e))$. Then $\xi:G\emb K$, and we write 
$\xi=\eta\circ\zeta$. 

If $G_0,G_1,\ldots,G_n$ are as in the paragraph before the previous one, 
then for each $i=1,2,\ldots,n$ there is a
natural homeomorphic embedding $G_{i-1}\emb G_i$,
and hence there is a
natural homeomorphic embedding $\iota:G_0\emb G_n$,
called the \dfn{canonical embedding determined by the
generating sequence $G_1,G_2,\ldots,G_n$}.
When there is no danger of confusion we will drop the reference to 
the generating sequence and simply talk about a canonical embedding.

\junk{
Let us assume the notation of (1.3). The result can be restated 
by saying that if there exists a homeomorphic embedding
$\eta:G\emb H$, then there exist a graph $G'$ as in (1.3) and
a homeomorphic embedding $\eta':G'\emb H$. We shall formulate
some of our results in a stronger form by insisting that $\eta'$ can
be chosen to be ``parallel" to $\eta$. 
This extra condition is often helpful in applications, because it allows
us to deduce that during the transition from $\eta$ to $\eta'$ certain
subgraphs stay fixed.
We now introduce this technical concept. 
}

We will frequently need to construct new \he s from old ones by means
of ``rerouting". We now introduce these constructions formally.
Let $G,H$ be graphs, and let $\eta:G\emb H$ be a 
homeomorphic embedding. Let $e\in E(G)$, and let $P'$ be a path
in $H$ of length at least one with both ends on $\eta(e)$, and otherwise disjoint
from $\eta(G)$. Let $P$ be the subpath of $\eta(e)$ with ends
the ends of $P'$. Let $\eta'(e)$ be the path obtained from
$\eta(e)$ by replacing the interior of $P$ by $P'$, and let $\eta'(x)=
\eta(x)$ for all $x\in V(G)\cup E(G)-\{e\}$. Then $\eta':G\emb H$
is a homeomorphic embedding, and we say that 
$\eta'$ was obtained from $\eta$
by \dfn{rerouting $\eta(e)$ along $P'$}.

Let $e,f,g$ be three distinct edges of $G$, all incident with
a vertex $v$ of degree three. Let $x$ be an interior vertex of  $\eta(e)$, let
$y$ be an interior vertex of $\eta(f)$, and let
$P'$ be a path in $H$ with ends $x$ and $y$, and otherwise disjoint
from $\eta(G)$. Let $\eta'(v)=y$, let $\eta'(e)$ be obtained
from $\eta(e)$ by deleting the part from $x$ to $\eta(v)$
(including $\eta(v)$ but not $x$) and adding $P'$, let
$\eta'(f)$ be obtained from $\eta(f)$ by deleting the part
from $y$ to $\eta(v)$ (including $\eta(v)$ but not $y$), and let
$\eta'(g)$ be obtained from $\eta(g)$ by adding the subpath of
$\eta(f)$ with ends $y$ and $\eta(v)$. For $z\in V(G)\cup E(G)
-\{v,e,f,g\}$ let $\eta'(z)=\eta(z)$. Then $\eta':G\emb H$, and
we say that 
$\eta'$
was obtained from $\eta$ by \dfn{rerouting $\eta(e)$ along $P'$}.

Let $e$ be an edge of $G$ with ends $u,v$ of degree three, 
let $f_1,f_2$ be the 
other two edges incident with $u$, and let $g_1,g_2$ be the other
two edges incident with $v$. 
Let $x$ be an interior vertex of  $\eta(f_1)$, let
$y$ be an interior vertex of $\eta(g_1)$, and let
$P'$ be a path in $H$ with
ends $x$ and $y$, and otherwise
disjoint from $\eta(G)$. Let $\eta'(u)=x$, let $\eta'(v)=y$, let
$\eta'(e)=P'$, let $\eta'(f_1)$ be the path obtained from 
$\eta(f_1)$ by deleting the subpath between $x$ and $\eta(u)$
(including $\eta(u)$ but not $x$), let $\eta'(g_1)$ be the 
path obtained from $\eta(g_1)$ by deleting the subpath between
$y$ and $\eta(v)$ (including $\eta(v)$ but not $y$), let
$\eta'(f_2)$ be obtained from $\eta(f_2)$ by adding the subpath
of $\eta(f_1)$ between $x$ and $\eta(u)$, and let $\eta'(g_2)$
be obtained from $\eta(g_2)$ by adding the subpath of $\eta(g_1)$
with ends $y$ and $\eta(v)$. For $z\in V(G)\cup E(G)-\{u,v,
e,f_1,f_2,g_1,g_2\}$ let $\eta'(z)=\eta(z)$. Then $\eta':G\emb H$,
and we say that 
$\eta'$ was obtained from $\eta$ by \dfn{rerouting 
$\eta(e)$ along $P'$}.

\junk{
Let $\eta,\eta':G\emb H$. We say that $\eta,\eta'$ are \dfn{close}
if they are $i$--close for some $i\in\{0,1,2\}$. We say that
$\eta,\eta'$ are \dfn{parallel} if for some integer $n>0$ there
exist homeomorphic embeddings $\eta_i:G\emb H$ $(i=1,2,\ldots,n)$
such that $\eta_1=\eta$, $\eta_n=\eta'$ and for $i=2,3,\ldots,n$,
$\eta_{i-1},\eta_i$ are close. Let $\eta:G\emb H$, let $\eta':
G'\emb H$, and let $G'$ be an extension of $G$. We say that 
$\eta,\eta'$ are \dfn{parallel} if $\eta$ is parallel to
the restriction of $\eta'$ to $G$.
}

\junk{
\newthm{3.1}Let $G,H,K$ be cubic graphs, and let $\eta_1,\eta_2:G\emb H$,
$\zeta_1,\zeta_2:H\emb K$ be homeomorphic embeddings. If $\eta_1$ is 
parallel to $\eta_2$ and $\zeta_1$ is parallel to $\zeta_2$, 
then $\eta_1\circ\zeta_1$ is parallel to $\eta_2\circ\zeta_2$.
In particular, if $H$ is an extension of $G$, then the restrictions
of $\zeta_1$ and $\zeta_2$ to $G$ are parallel.
\proof It is straightforward to verify that $\eta_1\circ\zeta_1$
is parallel to $\eta_1\circ\zeta_2$, and that  
$\eta_1\circ\zeta_2$ is parallel to $\eta_2\circ\zeta_2$, and hence
the first assertion follows. The second follows by taking
$\eta_1=\eta_2$ to be the canonical embedding. \qed
\newthm{3.2}Let $G,H$ be cubic graphs, 
let $u,v$ be two vertices of $G$ with neighbors $v,u_1,u_2$
and $u,v_1,v_2$, respectively,  for $i=1,2$ let
$G_i=G+(u,u_i,v,v_i)$, and let $k_i,l_i$ be the new vertices
of $G_i$. Then $G_1$ is isomorphic to $G_2$ under
an isomorphism $\phi:V(G_1)\to V(G_2)$ such that $\phi(u)=k_2$,
$\phi(k_1)=u$, $\phi(v)=l_2$, $\phi(l_1)=v$, and $\phi(x)=x$
for every $x\in V(G_1)-\{u,v,k_1,l_1\}$.
Moreover, if $\eta_1:G_1\emb H$ is a \he, then there exists
a \he\ $\eta_2:G_2\emb H$ such that the restrictions of $\eta_1$
and $\eta_2$ to $G$ are parallel, and $\eta_1(x)=\eta_2(\phi(x))$
for every $x\in V(G_1)$.
}

Our next objective is to analyze augmenting paths relative to \he s.
The next lemma follows by a standard application of network
flow theory. A proof may be found in [\cite{DieGT}, Lemma~3.3.3].

\newthm{3.3}Let $k\ge0$ be an integer, let $G,H$ be cubic graphs,
let $\delta_G A=\{e_1,e_2,\ldots,e_k\}$ be an edge-cut of
$G$ of cardinality $k$, and for $i=1,2,\ldots,k$ let the ends of $e_i$ be
$u_i\in A$ and $v_i\in V(G)-A$. 
Let $\eta:G\emb H$ be a 
homeomorphic embedding,
and assume that there is no edge-cut $\delta_H B$ of $H$ of cardinality
$k$ with $\eta(A)\subseteq B$ and $\eta(V(G)-A)\subseteq V(H)-B$.
Then there exist an integer $n$ and
disjoint paths $Q_1,Q_2,\ldots,Q_n$ in $H$, where $Q_i$ has distinct
ends $x_i$ and $y_i$ such that
\item{(i)}$x_1\in V(\eta(G|A))- \{\eta(u_1),\eta(u_2),\allowbreak\ldots,
\allowbreak\eta(u_k)\}$ and 
$y_n\in V(\eta(G\backslash A))-\{\eta(v_1),\allowbreak \eta(v_2),\allowbreak
\ldots,\eta(v_k)\}$,
\item{(ii)}for all integers $i\in\{1,2,\ldots,n-1\}$, the vertices $x_{i+1},y_{i}\in V(\eta(e_t))$
for some $t\in\{1,2,\ldots,k\}$, and $\eta(u_t),x_{i+1},\allowbreak y_i,
\eta(v_t)$ are
pairwise distinct and occur on $\eta(e_t)$ in the order listed, 
\item{(iii)}if $x_i,y_j\in V(\eta(e_t))$ for some
$t\in\{1,2,\ldots,k\}$
and $i,j\in\{1,2,\ldots,n\}$ with $i>j+1$, then
$\eta(u_t),y_j,x_i,\eta(v_t)$ occur on $\eta(e_t)$ in the order listed, and
\item{(iv)}for $i=1,2,\ldots,n$, if a vertex of $Q_i$ 
belongs to $V(\eta(G))$, 
then it is an end of $Q_i$.

In the situation described in~\refthm{3.3} we  call the sequence of paths
$\gamma=(Q_1,Q_2,\ldots,Q_n)$ an \dfn{augmenting sequence
with respect to $G,H,A$ and $\eta$}. 
Let $F$ be a subgraph of $G$.
We say that $\gamma$ is \dfn{reduced  modulo $F$} if the
following conditions are satisfied:
\item{(i)}If $e\in E(G|A)$ and $t\in\{1,2,\ldots,k\}$
are such that $x_1\in V(\eta(e))$ and $y_1\in V(\eta(e_t))$, then
$e$ and $e_t$ have no common end, and no end of $e$
is adjacent to an end of $e_t$ in $G\backslash E(F)$.
\item{(ii)}If $t\in\{1,2,\ldots,k\}$ and $f\in E(G\backslash A)$ 
are such that $x_n\in V(\eta(e_t))$ and $y_n\in V(\eta(f))$, 
then $e_t$ and $f$ have no common end, and no end of $e_t$
is adjacent to an end of $f$  in $G\backslash E(F)$.
\item{(iii)}If $t,t'\in\{1,2,\ldots,k\}$ and $i\in\{2,3,\ldots,n-1\}$ are
such that $x_i\in V(\eta(e_t))$ and $y_i\in V(\eta(e_{t'}))$, then
$t\ne t'$, $u_t$ is not adjacent to $u_{t'}$  in $G\backslash E(F)$, and $v_t$ is not adjacent
to $v_{t'}$  in $G\backslash E(F)$.

Let $G,H$ be graphs, let $\eta:G\emb H$, 
and let $F$ be a graph of minimum degree at least two (which includes the possibility 
that $F$ is empty).
We say that the \he\ $\eta$ \dfn{fixes} $F$ if 
$F$ is a subgraph of both $G$ and $H$,
$\eta(v)=v$ for every vertex
$v\in V(F)$ and for every edge $e\in E(F)$ the image $\eta(e) $ is the
path with edge-set $\{e\}$.
In many of our lemmas and theorems we will be able to find a \he\ that fixes
a specified graph $F$. This feature will not be needed in this
or the follow-up paper~[\cite{RobSeyThoCubic}], but is included because 
it may be useful in future applications.
As far as this paper and~[\cite{RobSeyThoCubic}] are concerned, the reader may take
$F$ to be the null graph.

The lemma we need is the following.

\newthm{3.4}Let $G,H$ be cubic graphs,
let $\delta_G A$
 be an edge-cut in $G$ such that no two members of $\delta_G A$ have a 
common end, let $F$ be a graph of minimum degree at least two, 
let  $\eta:G\emb H$ be a  homeomorphic embedding that fixes $F$,
and let $\gamma$ be an augmenting sequence with respect to
$G,H,A$ and $\eta$ of length $n$. 
Let us assume that $\gamma$ is minimal in the sense that 
there is no homeomorphic embedding $\eta':G\emb H$ that fixes $F$ and an
augmenting sequence 
with respect to
$G,H,A$ and $\eta'$ of length $n'$ such that $n'<n$.
Then $\gamma$ is reduced modulo~$F$.

\proof Let $G,H,A,\eta$ and $\gamma$  be as
stated, let $\delta_G A=\{e_1,e_2,\ldots,e_k\}$, 
and let $\gamma=(Q_1,Q_2,\allowbreak\ldots,Q_n)$. 
To prove that $\gamma$ satisfies (i) let $t$ and $e$ be as in (i), and
suppose for a contradiction that either $e$ and $e_t$ have a common
end, or that some end of $e$ is adjacent to some end of $e_t$  in $G\backslash E(F)$.
Since $G$ and $H$ are cubic,
$x_1$ is an interior vertex of
$\eta(e)$ and  $y_1$ is an interior vertex of $\eta(e_t)$.
Let $\eta'$ be obtained from
$\eta$ by rerouting $\eta(e)$ (if $e$ and $f$ have a common end)
or $\eta(g)$ (where $g$ is an edge of  $G\backslash E(F)$  adjacent to
both $e$ and $f$) along $Q_1$. 
Then $Q_2,Q_3,Q_4,\ldots,Q_n$ is an augmenting
sequence with respect to $G,H,A$ and $\eta'$, and hence $\gamma$ is not
minimal, a contradiction. 

Condition (ii) follows similarly, and so it remains to prove (iii).
To that end let $t,t'$ and $i$ be as in (iii). Suppose first that $t=t'$.
Then $\eta(u_t),x_i,y_{i-1},x_{i+1},y_i,\eta(v_t)$ all belong to $\eta(e_t)$ and
occur on $\eta(e_t)$ in the order listed.
Let $\eta'$ be obtained from $\eta$ by rerouting 
$\eta(e_t)$
along $Q_i$, and let $Q$ be the union of 
$Q_{i-1}$, $Q_{i+1}$ and the subpath of $\eta(e_t)$ with ends $y_{i-1}$
and $x_{i+1}$. Then $Q_1,Q_2,\ldots,Q_{i-2},Q,Q_{i+2},\ldots,Q_n$
is an augmenting sequence with respect to  $G,H,A$ and $\eta'$, and hence $\gamma$ is not minimal, a contradiction.

Thus $t\ne t'$. Next we suppose for a contradiction that $u_t$ is
adjacent to $u_{t'}$  in $G\backslash E(F)$. Let $\eta'$ be obtained from $\eta$
by rerouting $\eta(u_tu_{t'})$ along $Q_i$; then 
$Q_{i+1},Q_{i+2},\ldots,Q_n$
is an augmenting sequence with respect to  $G,H,A$ and $\eta'$, and hence $\gamma$ is not minimal, a contradiction. Similarly we deduce that
$v_t$ is not adjacent to $v_{t'}$. Thus $\gamma$ is reduced, as required.
\qed

Let $G,H$ be cubic graphs, 
let $\eta:G\emb H$ be a homeomorphic embedding,
let $e_1,e_2$ be two edges of $G$ with ends $u_1,v_1$ and $u_2,v_2$,
respectively, where $u_1,v_1,u_2,v_2$ are pairwise distinct,
and assume that there exists a path $Q$ in $H$ with ends $x_i
\in V(\eta(e_i))$ ($i=1,2$) and otherwise disjoint from $\eta(G)$.
Let
$G'=G+(u_1,v_1,u_2,v_2)$, and let $k_1,k_2$ be the new vertices of $G'$; 
then $G'$ is a $1$--extension of $G$. For $i=1,2$ let
$\eta'(k_i)=x_i$, let $\eta'(k_1k_2)=Q$,
let $\eta'(u_ik_i)$ be the subpath
of $\eta(u_iv_i)$ with ends $\eta(u_i),x_i$, let $\eta'(v_ik_i)$ 
be defined analogously, and let $\eta'(x)=\eta(x)$
for all $x\in V(G)\cup E(G)-\{e_1,e_2\}$.
Then $\eta':G'\emb H$ is a homeomorphic embedding. We say that
the pair $G',\eta'$ was obtained from $\eta$ by \dfn{routing the new 
edge along} $Q$. 


\newthm{3.5}Let $G$ be a cubic graph, let $H$ be a \cfc\ cubic graph,
let $F$ be a graph of minimum degree at least two,
let $\eta:G\emb H$ fix $F$, 
let $C$ be a quadrangle in
$G$ that is disjoint from $F$, and assume that $G$ has a circuit disjoint
from $C$. Then there exist a $1$--extension $G'$ of $G$  based at $C $
and a homeomorphic embedding $\eta':G'\emb H$ that fixes~$F$.

\proof
Since $H$ is \cfc\ and $G\backslash V(C)$ contains
a circuit, by \refthm{3.3} there exists an augmenting
sequence $\gamma=(Q_1,Q_2,\ldots,Q_n)$ with respect to $G,H,V(C)$ and $\eta$. 
By~\refthm{3.4} we may assume that  $\gamma$ is reduced modulo $F$.
Let $G',\eta'$ be obtained from $\eta$ by routing the new edge along $Q_1$.
Then  $G',\eta'$  satisfy the conclusion of \refthm{3.5}.~\qed

\newthm{3.6}Let $G,H$ be non-isomorphic cubic graphs, 
 let $F$ be a graph of minimum degree at least two,
let $\eta:G\emb H$ fix $F$, and let $G$ and $H$ be cyclically $4$--connected.
Then there exist a $1$--extension $G'$ of $G$
and a homeomorphic embedding $\eta':G'\emb H$ such that $\eta'$
fixes $F$.

\proof Since $G$ is not isomorphic to $H$, and $H$ is cyclically
$4$-connected, there exists a path $P$ in $H$ with at least
one edge, with both ends
on $\eta(G)$, and otherwise disjoint from $\eta(G)$. Let
$x_1\in \eta(e_1)$ and $x_2\in\eta(e_2)$ be the ends of $P$, where
$e_1,e_2\in E(G)$. Let $u,u_1$ be the ends of $e_1$, and let $v,u_2$ be
the ends of $e_2$. If $u,u_1,v,u_2$ are pairwise distinct, then
$G'=G+(u,u_1,v,u_2)$ and $\eta'$ obtained from $\eta$ by
routing the new edge along $P$ satisfy the conclusion of the lemma.
We may therefore assume that  say $u=v$. The case
when $u_1=u_2$ can be reduced to the case $u_1\not=u_2$ by
a similar, though easier argument, and is  omitted.
Thus we assume that $u=v$ and $u_1\not=u_2$. 

Let $G_1$ be obtained from $G$ by subdividing $e_1$ and $e_2$ and joining
the new vertices by an edge. Let $v_1$ and $v_2$ be the new vertices 
of $G_1$ numbered so that $v_1$ resulted by subdividing $e_1$. 
Let $\eta_1:G_1\emb H$ be obtained by routing the new edge
along $P$, and let $A=\{u,v_1,v_2\}$.
Since $H$ is cyclically $4$--connected,
there exists, by \refthm{3.3}, an augmenting sequence 
 with respect to $G_1,H,A$ and $\eta_1$. 
By \refthm{3.4} we may assume,
by replacing $\eta_1$ by a different \he\ if necessary, that
there exists a path $Q_1$ that is the first term of a reduced 
augmenting sequence modulo $F$  with respect to $G_1,H,A$ and $\eta_1$. 
Let $x\in V(\eta_1(G_1|A))$
and $y\in V(\eta_1(G_1\backslash A))$ be the ends of $Q_1$; 
let $f\in E(G\backslash A)$ be 
such that $y\in V(\eta_1(f))$. 
From the symmetry between $e_1$ and $e_2$ we may  assume that
$x$ belongs to $\eta(e_1)\cup\eta_1(v_1v_2)$.
Thus $P\cup Q_1$ has a subpath $R$ with one end in $\eta(e_1)$,
the other end $y$ and otherwise disjoint from $\eta(G)$.
If $f$ is not incident with $u_1$, then the graph and \he\
obtained from $G$ by routing the new edge along $R$ are as desired.
Thus we may assume that $f$  is incident with $u_1$.
Let $\eta':G\emb H$ be obtained from $G$ by rerouting $\eta(e_2)[u,v_2]$
along $P$; then the graph and \he\ obtained from $\eta'$ by routing the new edge 
along $Q_1$ are as desired.~\qed

\junk{
The last result of this section will be used in [\cite{RobSeyThoCubic}].
However, it is not needed in this paper. Let $G,H$ be cubic graphs,
and let $F^*$ be a subgraph of both $G$ and $H$. We say that a
\he\ $\eta:G\emb H$ \dfn{fixes} $F^*$ if $\eta(x)=x$ for every
$x\in V(F^*)\cup E(F^*)$.
degree at least two such that $F^*$ is a subgraph of both $G$ and $H$, and
let $F=G\backslash E(F^*)$. Let $\eta:G\emb H$ be a \he\ fixing $F^*$,
let $\zeta$ be the restriction of $\eta$ to $F$, and let $\zeta':F\emb H$
be parallel to $\zeta$. For $x\in V(F)\cup E(F)$ let $\eta'(x)=\zeta'(x)$,
and for $e\in E(F^*)$ let $\eta'(e)=e$. Then $\eta':G\emb H$ is a homeomorphic
embedding fixing $F^*$.
\noindent The proof is easy and is omitted.
}

\newsection{quad} FIXING A QUADRANGLE

Let $G$ be a \qc\ cubic graph, and let $C$ be a quadrangle in $G$.
In this section we study the following problem:
If $H$ is \cfc\ and topologically 
contains $G$, is there a \qc\ cubic graph $G'$ such that $G'$ is obtained from $G$ by one of
a set of well-defined operations, $G$ is topologically contained in $G'$, $G'$ is
topologically contained in
$H$ and has fewer quadrangles  than $G$?
The following simple result gives a preliminary answer.
Let us recall that extensions were defined at the beginning of
Section~\refsec{he}.

\junk{
Let $G_0$ be a \qc\ graph, and let $C_0$ be a quadrangle in $G$.
We say that $G_n$ is an \dfn{$n$-extension} of $G_0$ \dfn{based at} $C_0$ if
there exists a sequence $G_1,C_1,G_2, C_2,\dots, G_n$ such that for
$i=1,2,\dots, n$, the graph $G_i$ is a 1-extension of $G_{i-1}$ based at
$C_{i-1}$, and if $i<n$ then this 1-extension is short and $C_i$ is the
new quadrangle in $G_i$.  We say that $G_n$ is a \dfn{short}
\dfn{$n$-extension} of $G_0$ if $G_n$ is a short $1$-extension of $G_{n-1}$, 
and we say that it is a \dfn{long $n$-extension} otherwise.
}

\newthm{extex} 
Let $G$ be a quad-connected cubic
graph, let $H$ be a cyclically $5$-connected cubic graph,
 let $F$ be a graph of minimum degree at least two, and
let $\eta:G\emb H$ fix $F$.
Let $C$ be a 
quadrangle in $G$ that is disjoint from $F$.  Then there exist an integer $n\ge1$,
a long $n$-extension $G'$ of $G$ based
at $C$ 
and a \he\ $\eta':G'\emb H$ that fixes~$F$.

\proof Let $n$ be the maximum integer such that there exists an $n$-extension
$G'$ of $G$ based at $C$ and a \he\ $\eta':G'\emb H$ that fixes $F$.
This is well-defined, because $2k\le |V(J)|-|V(G)|$ for every $k$-extension $J$
of $G$.  We claim that $G'$ is long.  To prove the claim suppose to the 
contrary that it is short, and let $C'$ be the new quadrangle of $G'$.  
It follows that $C'$ is disjoint from $F$.
By \refthm{3.5} there exists a 1-extension $G''$ of $G'$ 
based at $C'$ and a
\he\ $\eta'':G''\emb H$  that fixes $F $.
Then $G''$ is an $(n+1)$-extension of $G$
based at $C$, contrary to the choice of~$n$.
This proves our claim that $G'$ is a long extension of $G$,
and hence the pair $G',\eta'$ satisfies the conclusion of the lemma.~\qed

In the rest of this section we strengthen~\refthm{extex} in two ways:
we give a bound on the minimum integer $n$ that satisfies the conclusion
of~\refthm{extex}, and we give an explicit list of long
extensions based at $C $ such that one of them is guaranteed to satisfy 
\refthm{extex}.
We now introduce these extensions.

\goodbreak\midinsert
\vskip5.5in
\includegraphics{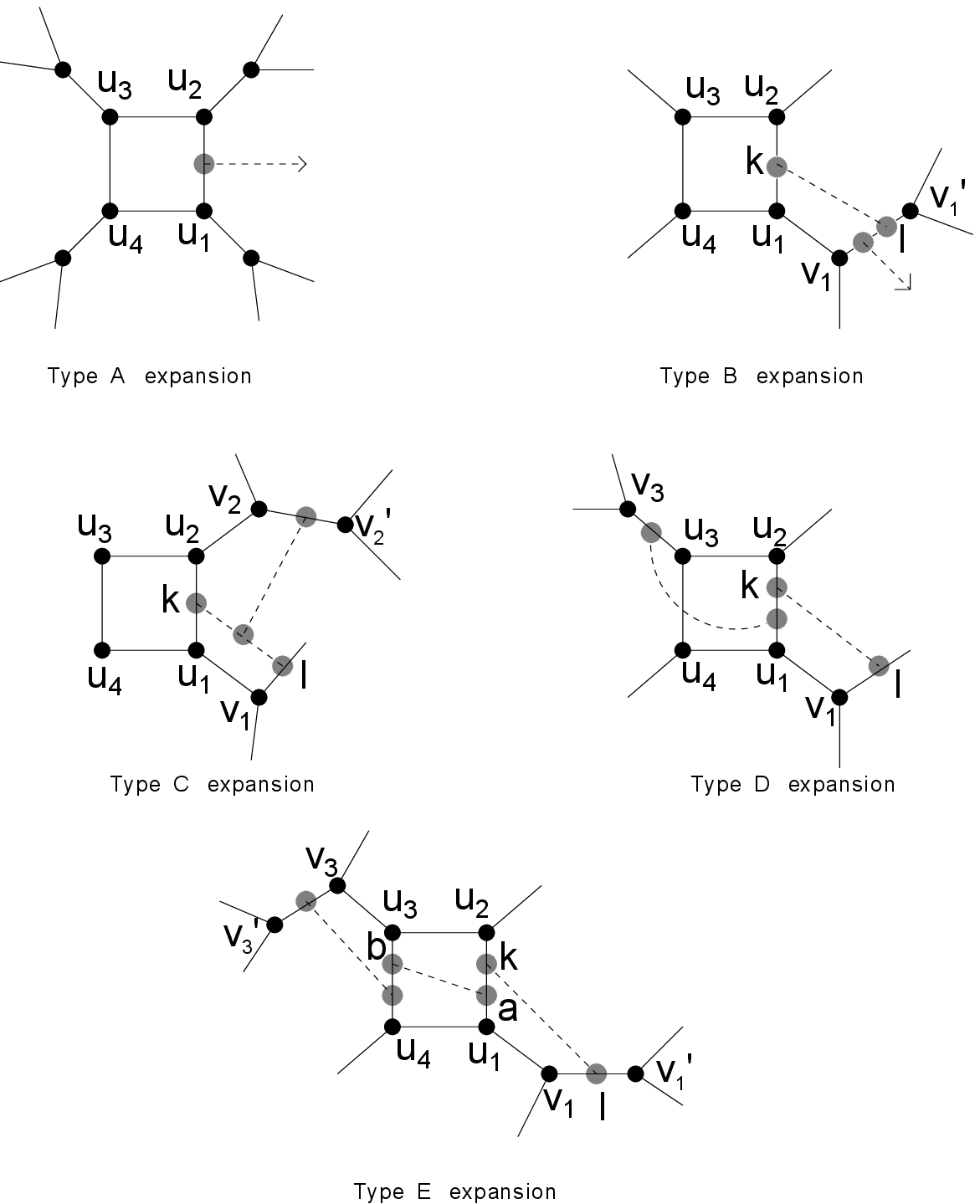}
\centerline{Figure 3. Type A--E expansions}
\endinsert

Let $G $ be a \qc\ cubic graph, let $C $ be a quadrangle  in $G$, let
 $u_1,u_2,u_3,u_4$ be the vertices of $C$ in order, for 
$i=1,2,3,4$ let $v_i$ be the unique neighbor of $u_i$ not on $C$,
and let $v_i'\not=u_i$ be a neighbor of $v_i$. 
It follows that $v_i'\not\in\{v_1,v_2,v_3,v_4\}$.
Let $G_1=G+(u_1,u_2,x,y)$ be a $1$--extension of $G$ with $x,y\not\in V(C)$. 

\item{$\bullet$}
If $G_1$ is a long extension of $G$,
we say that $G_1$ is a \dfn{type  A
expansion of $G$} based at $C$. 

\noindent
See Figures~3, 4 and 5. 
Otherwise we may assume  that say $x=v_1$ and $y=v_1'$.
Let $C_1$ be 
the new quadrangle  of $G_1$; 
thus $C_1$ has vertex-set $\{v_1,u_1,k,l\}$, where $k,l$
are the new vertices of $G_1$. 

\item{$\bullet$}
Let $G_2=G_1+(v_1,l,a,b)$ be
a $1$--extension of $G_1$. If $G_2$ is a long extension of $G$
we say that $G_2$ is a \dfn{type B expansion of $G$} based at $C$,
and that the sequence $G_1,G_2$ is a \dfn{standard generating sequence of $G_2$}.

\item{$\bullet$}
Let $G_3$ be  $G_1+(k,l,v_2,v_2')$ or $G_1+(u_1,v_1,v_4,v_4')$.
If $G_3$ is a long extension of $G$ we say that $G_3$ is a 
\dfn{type C expansion of $G$} based at $C$,
and that the sequence $G_1,G_3$ is a \dfn{standard generating sequence of $G_3$}.

\noindent
(We apologize for the double use of the letter C and hope it causes no confusion.)

\item{$\bullet$}
Let $G_4$ be the graph $G_1+(u_1,k,u_3,v_3)$; we say that $G_4$ is a
\dfn{type D expansion of $G$} based at $C$,
and that the sequence $G_1,G_4$ is a \dfn{standard generating sequence of $G_4$}.

\item{$\bullet$}
Let $G_5'$ be the graph $G_1+(u_1,k,u_3,u_4)$,
let $a,b$ be the new vertices of $G_5'$, and let $G_5$ be the
graph $G_5'+(b,u_4,v_3,v_3')$. The graph $G_5$ is called a \dfn{type
E expansion of $G$} based at $C$.


\goodbreak\midinsert
\vskip3in
\includegraphics{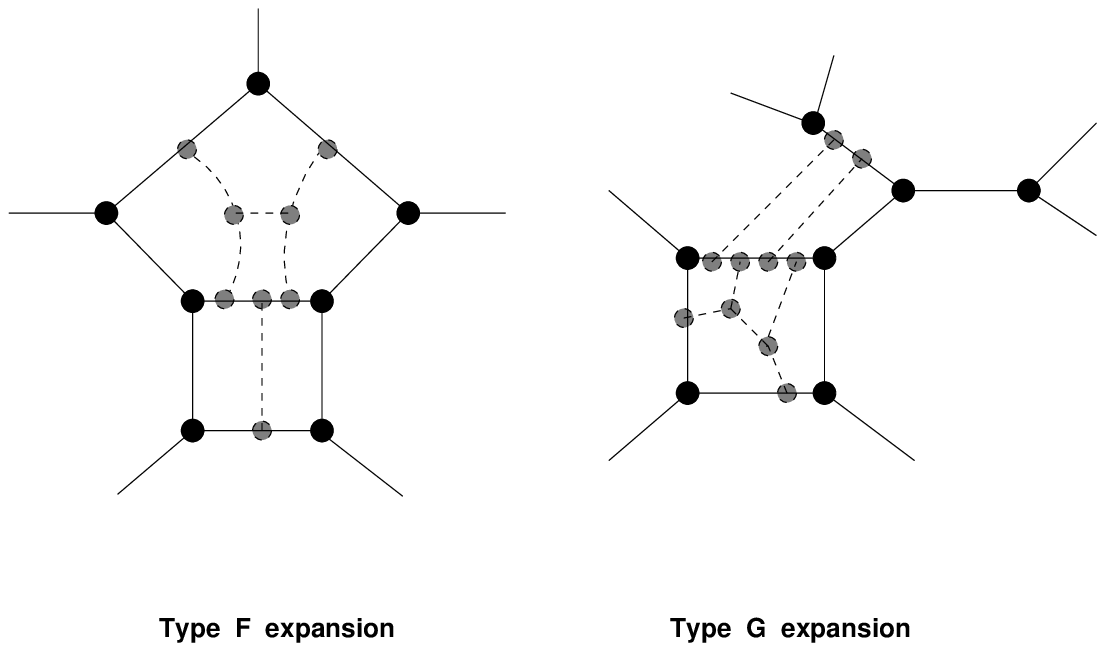}
\centerline{Figure 4. Type F and G expansions}
\endinsert

\item{$\bullet$}
If $v_1'=v_2'$, then 
let $G_6'=G_1+(u_2,k,v_2,v_2')$, let $k_2,l_2$ be the new vertices
of $G_6'$, let $G_6''=G_6'+(k_2,k,u_3,u_4)$, 
and let $G_6=G_6''+(k,l,k_2,l_2)$. 
The graph $G_6$ is called a \dfn{type
F expansion of $G$} based at $C$. We also say that $G_6$
is a \dfn{type F expansion of $G$ based on}
 $(u_1,u_2)$,
and that $k,l,k_2,l_2,k_3,l_3,k_4,l_4$ (in the order listed)
are the \dfn{new vertices} of $G_6$, where $k_3,l_3$ are the new
vertices of $G_6''$ and $k_4,l_4$ are the new vertices of $G_6$.
We say that $$\{u_1,u_2,u_3,u_4,v_1,v_2,v_1',k,l,k_2,l_2,k_3,l_3,k_4,
l_4\}$$ is the \dfn{core} of the type F expansion $G_6$.

\item{$\bullet$}
Let $G_7$ be a type F expansion of $G_1$ 
based on $(u_1,k)$. We say that $G_7$ is a \dfn{type G expansion
of $G$} based at $C$. 

\noindent
(Again, apologies for the double use of the letter G.)

\item{$\bullet$}
Assume now that $G$ has a quadrangle $D$ with vertex-set $x_1,x_2,x_3,x_4
\in V(G)-V(C)$ such that $x_1$ is adjacent to $u_1$, 
the vertices $u_2$ and $x_2$ have a 
common neighbor, and $u_4$ and $x_4$ have a common neighbor.
Assume further that $\{x,y\}=\{x_1,x_2\}$, and let us recall that $k,l$ are
the new vertices of $G_1$.  Let $G_8$ be a type~F expansion of $G_1$ based 
on $(x_1,l)$; in those circumstances we say that $G_8$ is a \dfn{type~H
expansion of $G$}.  
The \qcity\ of $G $ implies that in this case $|V(G)|\le14$.

\goodbreak\midinsert
\vskip3in
\includegraphics{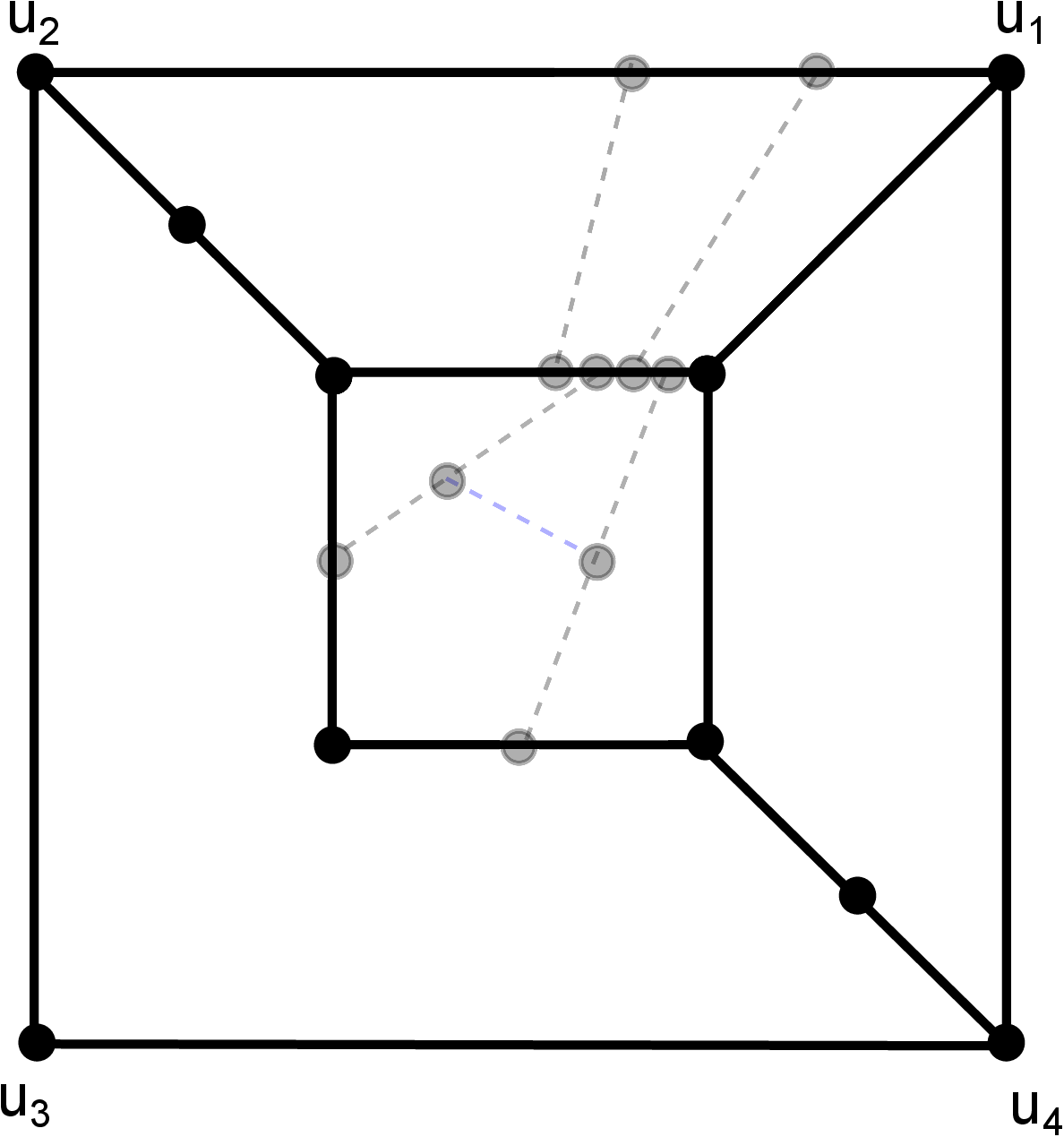}
\centerline{Figure 5. Type H expansion}
\endinsert

\noindent
It follows from the \qcity\ of $G$ that
$G_4$, $G_5$, $G_6$, $G_7$, $G_8$ are long extensions of~$G$.  
We offer the following easy but important remark.
Let us recall that generating sequences were defined at the beginning of Section~\refsec{he}.



\newthm{ezremark}
Let $G_0$ be a \qc\ cubic graph, let $C_0$ be a quadrangle in $G_0$, let
$G'$ be a type B, C, D, E, F, G or H expansion of $G_0$ based at $C_0$,
let $G_1,G_2,\ldots,G_k$ be a generating sequence of $G'$ from $G_0$
based at $C_0$, and let $F$ be a graph of minimum degree at least two
disjoint from $C_0$.
Let the vertices of $C_0$ be $u_1,u_2,u_3,u_4$ in order, 
let  $v_1$ be the neighbor of $u_1$ not on $C_0$,
and let $v_1', v_1'' $ be the two neighbors of $v_1$ other than $u_1$.
If $G_1=G_0+(u_1,u_2,v_1,v_1')$, then there exists a generating sequence
$G_1',G_2',\ldots,G_k'$ of $G'$ from $G_0$ based at $C_0$ such that
\item{$\bullet$} for $i=1,2,\ldots,k$ the graph $G_i'$ is isomorphic to $G_i$,
\item{$\bullet$}
$G_1'=G_0+(u_1,u_4,v_1,v_1'')$,
\item{$\bullet$}
if $F$ is a subgraph of both $G_0 $ and $G_k $, then $F$ is a subgraph of
$G_k' $,
and
\item{$\bullet$}
if  the sequence $G_1,G_{2}$
is a standard generating sequence of $G_2$, then
the sequence $G_1',G_{2}'$
is a standard generating sequence of~$G_2'$.

\noindent
The proof is clear.





\newthm{4.1}Let $F$ be a graph of minimum degree at least two, let
$G$ be a \qc\ cubic graph, let $C$ be
a quadrangle in $G$, and let $G_2$ be a long $2$-extension of $G$ based at $C$
such that $F $ is a subgraph of both $G $ and $G_2$
and $F$ is disjoint from $C$.
Then there exist an expansion $G'$ of $G$ of type A, B, C, or D based at $C$,
and a \he\ $\eta':G'\emb G_2$ such that $\eta'$ fixes~$F $.

\proof Let $G,C,G_2$ be as stated, let $u_1,u_2,u_3,u_4$ be the
vertices of $C$ (in order), and let $v_i$ be the neighbor of $u_i$ not on $C$.
Let $G_2=G_1+(u,v,x,y)$ and
$G_1=G+(u_1,u_4,v_1,v'_1)$, where 
 $v'_1\not\in V(C)$ is adjacent to $v_1$,
and $\{u,v\}$ is one of $\{v_1,l_1\}$, $\{k_1,l_1\}$, $\{u_1,v_1\}$,
$\{u_1,k_1\}$, where $k_1,l_1$ are the new vertices of $G_1$.
Let $k_2,l_2$ be the new vertices of $G_2$.

First, if $\{u,v\}=\{v_1,l_1\}$, then $G_2$ is a type B expansion
of $G$, and hence $G_2$ and the identity \he\ $G_2\emb G_2$
satisfy the conclusion
of the lemma. Second, let us assume that $\{u,v\}=\{k_1,l_1\}$.
By considering the path $k_1k_2l_2$ we see that there exists a \he\ 
$G+(u_1,u_4,x,y)\emb G_2$ that fixes $F$, and hence we may assume that the $1$-extension
$G+(u_1,u_4,x,y)$ is short.
It follows that 
$\{x,y\}$ equals one of $\{v_4,v_4'\}$, $\{u_3,v_3\}$,
$\{u_2,u_3\}$ or $\{u_2,v_2\}$, where $v_4'\not=u_4$ is a neighbor of
$v_4$. We break the analysis into three subcases. First, if
$\{x,y\}=\{v_4,v_4'\}$, then $G_2$ is a type C expansion of $G$,
and hence $G_2$ and the identity \he\ $G_2\emb G_2$
satisfy the conclusion
of the lemma.
For the second subcase assume that
$\{x,y\}=\{u_3,v_3\}$ or $\{u_2,u_3\}$. Let 
$G'=G+(u_3,u_4,v_1,v_1')$; then $G'$ is a long $1$-extension of $G$ by 
\refthm{2.6}, and hence
is a type A expansion of $G$. 
Let $\eta:G\emb G_2$ be the canonical \he\ determined by the generating sequence
$G_1, G_2$.
Let $\eta':G'\emb G_2$ be obtained from $\eta$ first by rerouting
$u_3u_4$ along $k_1k_2l_2$, and then
routing the new edge along $k_2l_1$. 
Since $u_1,u_2,u_3,u_4,k_1,l_1,k_2,l_2\not\in F$ we deduce that $\eta'$ fixes $F$. 
(In the future we will omit this kind of argument,
because it will be clear that all the \he s that we will construct will fix $F$.)
The pair $G',\eta'$
satisfies the conclusion of the lemma.
The third and last subcase
is that $\{x,y\}=\{u_2,v_2\}$. Let $G'=G+(u_2,u_3,v_1,v_1'')$,
where $v_1''\not\in\{v_1',u_1\}$ is the third neighbor of $v_1$.
Then $G'$ is a long $1$-extension of $G$ by \refthm{2.6}. 
Let $\eta':G'\emb G_2$ be obtained
from $\eta$ first by rerouting $\eta(u_1v_1)$ along 
$k_1k_2l_1$,
then rerouting $\eta(k_1u_1u_2)$ along $k_2l_2$,
and finally routing the new edge along $\eta(u_2u_1v_1)$.
Then $G',\eta'$ satisfy the conclusion of the lemma.
This completes the case $\{u,v\}=\{k_1,l_1\}$.

The third case $\{u,v\}=\{u_1,v_1\}$ is symmetric to the previous 
case by \refthm{ezremark}, and so we proceed to the fourth and last case,
namely $\{u,v\}=\{u_1,k_1\}$. Let $G'=G+(u_1,u_4,x,y)$ and $\eta':G'\emb G_2$
be obtained from $G,\eta$ by routing the new edge along $\eta(k_2l_2)$.
We may assume that $G'$ is a short $1$-extension of $G$, for otherwise
the lemma holds. Thus  either $\{x,y\}=\{u_3,v_3\}$,
or $\{x,y\}=\{v_4,v_4'\}$, where $v_4'\not=u_4$ is a neighbor of 
$v_4$. In the former case $G_2$ is a type D expansion of $G$,
and so the lemma holds, and hence we may assume that the latter
case holds. Since by \refthm{ezremark} there is symmetry between $u_1v_1$
and $k_1l_1$ we deduce that also $\{x,y\}=\{v_2,v_2'\}$, where
$v_2'\not=u_2$ is a neighbor of $v_2$. It follows that $v_2$
and $v_4$ are adjacent in $G$, contrary to the \qcity\ of $G$.
This completes the fourth case, and hence the proof of the lemma.~\qed

\junk{
\noindent {\bf (4.2+1/2)}
{\sl
Let $G$ be a \qc\ graph, let $C$ be a quadrangle in $G$, and let 
$G_2$ be a type B expansion of $G$ based at $C$, 
where the notation is as in the definition of type B 
expansion. So, in particular, $G_1=G+(u_1,u_2,v_1,v_1')$ and
$G_2=G_1+(v_1,l,a,b)$.
Let $G_1'=G+(u_1,u_4,v_1,v_1'')$, where $v_1''$ is the third neighbor of $v_1$,
let $k',l'$ be the new vertices of $G_1'$, and let $G_2'=G_1'+(v_1,l',a,b)$.
Then $G_2'$ is isomorphic to $G_2$. 
}
}

\newthm{4.2}Let $F$ be a graph of minimum degree at least two, let
$G$ be a \qc\ cubic graph, let $C$ be
a quadrangle in $G$, and let $G_3$ be a long $3$-extension of $G$
based at $C$ such that $F$ is a subgraph of both $G $ and $G_3$
and $F$ is disjoint from $C$.
Then there exist a graph $G'$ and a \he\ $\eta':G'\emb G_3$ such that
$\eta'$ fixes $F$ and $G'$ is either a type E expansion
or a long $1$- or $2$-extension of $G$ based at $C$.


\proof
Let $G_1$ be a short $1$-extension of $G$ based at $C$
such that $G_3$ is a long $2$-extension of $G_1$ based at the new quadrangle 
$C_1$ of $G_1$.
By \refthm{4.1} applied to the graph $G_1$ and circuit $C_1$ there exist 
an expansion $G_3'$ of $G_1$ of type A, B, C, or D based at $C_1$,
and a \he\ $\eta_3':G_3'\emb G_3$ that fixes $F$.
If $G_3'$ is of type A, then it is a long $2$-extension of $G$ based at $C$,
and the lemma holds. 
Thus we may assume that $G_3'$ is of type B, C, or D.
It follows that $G_3$ and $G_3'$ have the same number of vertices,
and hence $\eta_3'$ is an isomorphism of $G_3$ and $G_3'$.
It follows that if the conclusion of the lemma holds for $G_3'$,
then it holds for $G_3$. Therefore
we may assume that $G_3'=G_3$ and that $\eta_3'$
is the identity \he.
In other words, $G_3$ is a  type B, C, or D expansion of $G_1$ based at $C_1$.

Let $G_2,G_3$ be a standard generating sequence of the expansion $G_3$.
Then $G_1,G_2,G_3$ are \qc\ by \refthm{1extquadconn}.
Let the vertices of $C$ be $u_1,u_2,u_3,u_4$ in order. 
For $i=1,2,3,4$ let $v_i$ be the neighbor of $u_i$ not on $C$,
let $v_i',v_i''$ be the neighbors of $v_i$ different from $u_i$,
and let $w_i,w_i'$ be the neighbors of $v_i''$ different from $v_i$.
Let $G_2=G_1+(a_1,a_2,a_3,a_4)$. 
We may assume that $G_1=G+(u_1,u_4,v_1,v'_1)$.
Let $u_5,v_5$ be the new vertices of $G_1$; then
$V(C_1)=\{v_1,u_1,u_5,v_5\}$.
We claim the following.

\claim{1}We may assume that $\{a_1,a_2\}=\{u_1,v_1\}$, and that
$\{a_3,a_4\}$ is equal to one of $\{u_2,u_3\}$, $\{u_2,v_2\}$ or
$\{v_1'',w_1'\}$.

To prove (1) we first note that by \refthm{ezremark} there is  symmetry
between $\{u_5,v_5\}$ and $\{u_1,v_1\}$, and so we may assume
that $\{a_1,a_2\}\not=\{u_5,v_5\}$.
Secondly, assume that
$\{a_1,a_2\}=\{v_1,v_5\}$. Since $G_2$ is not a long extension of $G_1$, one
of $a_3,a_4$ equals one of $v_1',v_1''$. Let us assume that
$a_3=v_1''$; the argument for $v_1'$ is symmetric. We may assume
from the symmetry that $a_4=w_1'$. 
It follows from \refthm{ezremark} applied to the graph $G_1$ and cycle $C_1$
that we may replace $G_2$ by the graph
$G_1+(u_1,v_1,v_1'',w_1)$ and thus arrange for
the first assertion of (1) to hold. The case 
$\{a_1,a_2\}=\{u_1,u_5\}$ follows similarly. This proves that we may
assume that $\{a_1,a_2\}=\{u_1,v_1\}$. 
Since $G_2$ is a short extension of $G_1$, 
we see that $\{a_3,a_4\}$ is equal to one of $\{u_2,u_3\}$, 
$\{u_2,v_2\}$, $\{v_1'',w_1'\}$ or $\{v_1'',w_1\}$. Since the last
two cases are symmetric, we may assume that one of the first
three occurs. This proves (1).
\medskip

Let $k_2,l_2$ be the new vertices of $G_2$, let 
$G_3=G_2+(a_5,a_6,x,y)$, and let $k_3,l_3$ be the new vertices of
$G_3$.
Let $\eta:G\emb G_3$ be the canonical \he\ determined by the generating sequence
$G_1, G_2,G_3$.
Since $F$ has minimum degree at least two and  is a subgraph of both $G $ and $G_3$
 we deduce that

\claim{2}$u_1,u_2,u_3,u_4,u_5,v_1,v_5,k_2,l_2,k_3,l_3\not\in F$.

To make the forthcoming case analysis 
easier to follow let us make an outline. There will be three
supercases depending on $\{a_3,a_4\}$. These will be divided
into cases depending on the type of the expansion $G_3$, and
the cases will sometimes be further divided into subcases
depending on $G_3$. In each subcase we shall construct a pair
$G',\eta'$ that satisfies the conclusion of the theorem.
We first dispose of the supercase
$\{a_3,a_4\}=\{u_2,u_3\}$. Let $G'=G+(u_2,u_3,v_1,v_1'')$, and let
$\eta':G'\emb G_3$ be obtained from $\eta$ by rerouting $u_1v_1$ along $u_5v_5$ 
and then  routing the new edge along  the path $l_2k_2v_1$. 
By \refthm{2.6} and (2) the pair $G',\eta'$ satisfies the conclusion of the lemma.
This completes the first supercase.

For the second supercase we assume that $\{a_3,a_4\}=\{u_2,v_2\}$.
This will be divided into cases. As a first case assume
that $G_3$ is a type B expansion of $G_1$. Then
$\{a_5,a_6\}=\{l_2,u_2\}$. Assume as a first subcase that
$\{x,y\}$ is not equal to any of $\{v_3,v_3'\}$, $\{v_3,v_3''\}$,
$\{u_4,v_4\}$, or $\{u_4,u_5\}$.
%
Let $G',\eta'$ be obtained from $\eta$ 
first by rerouting $u_1u_2$ along $l_2k_2$ and then
by routing the new edge along $k_3l_3$ if $\{x,y\}\ne\{u_5,v_5\}$
and along $k_3l_3v_5$ otherwise. 
Then $G'$ is a long 1-extension of $G$, and so the lemma holds.
We may therefore assume that $\{x,y\}$ is equal
to one of the sets specified above. As a second subcase assume
that $\{x,y\}=\{v_3,v_3'\}$ or  $\{x,y\}=\{v_3,v_3''\}$. Then
$G_3$ is isomorphic to a type E expansion of $G$, and so the
lemma holds. Thirdly, let us assume that $\{x,y\}=\{u_4,v_4\}$ or
$\{x,y\}=\{u_4,u_5\}$.  Let 
$G'=G_1+(u_1,u_5,u_3,v_3)$; then $G'$ is a long 2-extension of $G$.
Thus $G'$ and the homeomorphic embedding
obtained from the canonical \he\ $G_1\emb G_3$ 
(determined by the generating sequence $G_2,G_3$) first
by rerouting $u_1u_2$ along $l_2k_2$, then
rerouting $u_3u_4$ along $k_3l_3$, and finally 
routing the new edge along $u_1u_2$ 
satisfy the conclusion of the lemma.
This completes the case when $G_3$ is a type B expansion.

\junk{
As a fourth subcase assume that 
$\{x,y\}=\{u_5,v_5\}$. Let $G'=G+(u_2,u_3,v_1,v_1')$; then
$G'$ is a long extension of $G$
by \refthm{2.6}, and hence $G'$ and the
homeomorphic embedding $G'\emb H$ obtained from $\zeta_2$ by 
routing the new edge along $\eta_3(k_3l_3)\cup\eta_3(l_3v_5)$ 
again contradict the $C$-complexity of $\eta$ being three.
The fifth and final subcase is that $\{x,y\}=\{k_2,v_1\}$. Let
$G'=G+(u_2,u_3,v_1,v_1'')$. Then $G'$ is a long extension of $G$
by \refthm{2.6}, and $G'$ and the homeomorphic embedding $G'\emb H$
obtained from $\zeta_1$ by rerouting $\eta_0(u_1u_2)$ along
$\eta_2(l_2k_2)\cup\eta_2(k_2u_1)$ and routing the new edge along
$\eta_3(k_3l_3)\cup\eta_3(l_3v_1)$ again result in a contradiction.
This completes the case when $G_3$ is a type B expansion.
}

For the second case assume that $G_3$ is a type C expansion of $G_1$.
There are two subcases. Assume first that $\{a_5,a_6\}=\{u_1,u_2\}$.
Then $\{x,y\}=\{u_4,v_4\}$, because $\{x,y\}\not=\{u_3,u_4\}$, since
$G_3$ is a long 1-extension. Let $G'=G+(u_3,u_4,v_1,v_1')$ and let $\eta':G'\emb G_3$
be obtained from $\eta$ first by rerouting $u_1u_4$
along $k_3l_3$, and then by routing the new edge along
$u_4u_5v_5$. The graph $G'$ is a long $1$-extension of $G$
 by \refthm{2.6}, and hence~\refthm{4.2} holds.
The second subcase is that $\{a_5,a_6\}=\{k_2,l_2\}$. Then say
$x=v_1''$ and $y$ is a neighbor of $v_1''$ different from $v_1$.
Since $G_3$ is a long extension of $G$ we deduce that $v_2\not=y$. Let $G'=
G+(u_1,u_2,v_1'',y)$ and let $\eta':G'\emb G_3$ be obtained from
$\eta$ first by rerouting $u_1u_2$ along $l_2k_2$ and then
routing the new edge along $k_3l_3$. Then
$G'$ is a long 1-extension of $G$ by \refthm{2.6},  and hence~\refthm{4.2} holds.
This completes the second case. For
the third case we assume that $G_3$ is a type D expansion of $G_1$.
Then $\{a_5,a_6\}=\{u_1,k_2\}$ and $\{x,y\}=\{v_1',v_5\}$.
Let $G'=G_1+(u_4,v_4,v_1,v_5)$; then $G'$ is a long 1-extension of $G_1$ by \refthm{2.6},
and hence it is a long $2$-extension of $G$.
Let $\eta'$ be obtained from the canonical \he\ $G_1\emb G_3$
first by rerouting $u_1u_2$ along $k_2l_2$,
then rerouting $u_5v_5$ along $k_3l_3$, then
rerouting $u_3u_4$ along $u_1u_2$, and finally
routing the new edge along $u_5v_5$. 
Again, the pair $G',\eta'$ satisfies the conclusion of the lemma.
This completes the third case and hence the second supercase.

The third and last supercase is that $\{a_3,a_4\}=\{v_1'',w_1'\}$.
We claim that we may assume that $w_1'=v_2$. Indeed, suppose that
$w_1'\not=v_2$, let $G'=G+(u_1,u_2,v_1'',w_1')$ and let $\eta'$
be obtained from $\eta$  by rerouting $u_1v_1$
along $u_5v_5$ and
by routing the new edge along
$u_1k_2l_2$. Then $G',\eta'$ 
satisfy the conclusion of the lemma by \refthm{2.6}.
This proves that we may assume that  $w_1'=v_2$. From the
symmetry we may assume that $v_2'=v_1''$. 
We distinguish three cases depending on whether $G_3$ is of type B, C, or D.

For the first case assume that $G_3$ is of type B.
Then $\{a_5,a_6\}=\{l_2,v_1''\}$.
Let us first dispose of the case when one of $x,y$ is equal to $v_2''$;
say $x=v_2''$. Then $y\ne v_2$, because $G_3$ is a long extension of $G $.
Let $G'$ and $\eta':G'\emb G_3$
be obtained from $\eta$ by first  rerouting $v_2v_2''$ along $k_3l_3$, 
then rerouting $u_1u_2$ along $k_2l_2$,
then routing the first new edge along $u_5v_5$,
then routing the second new edge along $u_1u_2$,
and finally routing the third new edge along $v_2v_2''$.
Then $G'$ is a type E expansion of $G$, 
and thus $G',\eta'$ satisfy the conclusion of
the lemma. 
We may therefore assume that $x,y\ne v_2'$.
Let us assume next that $\{x,y\}$ is not equal to any of the pairs
$\{u_1,u_5\}, \{u_1,u_2\}, \{u_2,u_3\}$.
Let $\xi$ be obtained from $\eta$ by first  rerouting $\eta(u_1v_1)$ along $u_5v_5$, 
and then rerouting $\eta(v_2v_2')$ along $l_2k_2v_1$.
Let $G''$ and $\eta'':G''\emb G_3$
be obtained from $\xi$ by routing the first new edge along $u_1k_2$,
and then routing the second new edge along $l_2v_2'$.
Then $G''$ is a long $2$-extension of $G$,  and hence the pair 
$G'',\eta''$ satisfies the conclusion of the lemma.
We may therefore assume that $\{x,y\}$ is equal to one of the pairs
$\{u_1,u_5\}, \{u_1,u_2\}, \{u_2,u_3\}$.
Let $G'''$ and $\eta''':G'''\emb G_3$
be obtained from $\xi$ by routing the new edge along $v_2'k_3l_3$.
Then $G'''$ is a long $1$-extension of $G$,  and hence the pair 
$G''',\eta'''$ satisfies the conclusion of the lemma.
This completes the first case.

For the second case assume that $G_3$ is of type C.
Then $\{a_5,a_6\}=\{v_1,v_1''\}$ or $\{a_5,a_6\}=\{k_2,l_2\}$.
Thus we distinguish two subcases.
Assume as a first subcase that $\{a_5,a_6\}=\{v_1,v_1''\}$. 
It follows that one of $x,y$ equals
$v_1'$, say $x=v_1'$. Then $y\not=v_5$; let $y'\not\in\{v_5,y\}$
be the third neighbor of $v_1'$. Let $\zeta'$ be obtained from
$\eta$ by rerouting $v_1v_1'$
along $k_3l_3$.
Let $G'=G+(u_1,u_4,v_1',y')$ and $\eta':G'\emb G_3$
be obtained from $\zeta'$ by routing the new edge along
$u_5v_5v_1'$; if $v_4\not=y'$ then $G',\eta'$
satisfy the lemma by \refthm{2.6}.  Thus we may assume that $y'=v_4$.
Let $G''=G+(u_1,u_4,v_1',v_4)+(u_2,u_3,v_2,v_2'')$ and let
$\eta''$ be obtained from $\zeta'$ first by rerouting
$u_1u_5$ along $v_1v_5u_5$, 
then by rerouting
$u_2v_2$ along $k_2l_2$, then routing the
first new edge along $v_5v_1'$, and then routing
the second new edge along $u_2v_2$.
Let $G',\eta'$ be obtained from $\eta''$ by routing the new edge along  $\eta_3(u_1u_5)$. 
Then $G'$ is a type E expansion of $G$, 
and thus $G',\eta'$ satisfy the conclusion of
the lemma. This completes the first subcase.
For the second subcase assume that $\{a_5,a_6\}=\{k_2,l_2\}$; then
$\{x,y\}=\{u_2,u_3\}$, because
$\{x,y\}\not=\{u_2,v_2\}$ by the fact that $G_3$ is a long 1-extension of $G$.
Let $G'=G+(u_2,u_3,v_1,v_1'')$ and let $\eta'$ be obtained
from $\eta$ by rerouting $u_1v_1$ along $u_5v_5$ and then routing the new edge along
$l_3k_3k_2v_1$. Then $G',\eta'$
satisfy the conclusion of the lemma by \refthm{2.6}.
This completes the second subcase and hence the second case.

For the third case assume that $G_3$ is of type D.
Then $\{a_5,a_6\}=\{v_1,k_2\}$ and $\{x,y\}=\{u_4,u_5\}$.
Let $G'=G+(u_1,u_2,v_2,v_2')+(a,u_2,v_1,v_1')$ 
(where $a,b$ are the new vertices of $G+(u_1,u_2,v_2,v_2')$)
and $\eta':G'\emb G_3$
be obtained from $\eta$ by  rerouting $\eta(u_1u_4)$ along
$k_3l_3$, then routing the first new edge along $k_2l_2$,
and then routing the second new edge along $u_1u_5v_5$.
Then $G'$ is a long $2$-extension of $G$ by \refthm{2.6},
and hence the pair $G',\eta'$ satisfies the conclusion of the lemma.
This completes the third case, and hence the third supercase,
and thus the proof of the lemma.~\qed

\newthm{4.3}Let $F$ be a graph of minimum degree at least two, let
 $G$ be a \qc\ cubic graph, let $C$ be
a quadrangle in $G$, and let $G_4$ be a long $4$-extension of $G$
based at $C$ such that $F$ is a subgraph of both $G $ and $G_4$
and $F$ is disjoint from $C$.
Then there exist a graph $G'$ and a \he\ $\eta':G'\emb G_4$ such that
$\eta'$ fixes $F$ and $G'$ is either a type F expansion of $G$ based at $C$
or a long $n$-extension of $G$ based at $C$ for some $n\in\{1,2,3\}$.

\proof 
Similarly as in the proof of \refthm{4.2} we may assume that there exists
a short 1-extension $G_1=G+(u_1,u_4,v_1,v'_1)$ of $G$ based at $C$ such 
that $G_4$ is a type E expansion of $G_1$ based at the new quadrangle $C_1$ of $G_1$.
%
Since $G_1$ is short one of $u_1,u_4$ is adjacent to
one of $v_1,v_1'$, and so we may assume that $u_1$ is adjacent
to $v_1$. Let the vertices of $C$ be $u_1,u_2,u_3,u_4$ in order,
and for $i=1,2,3,4$ let $v_i$ be the neighbor of $u_i$ not on $C$.
Let $u_5,v_5$ be the new vertices of $G_1$; then $V(C_1)=
\{v_1,u_1,u_5,v_5\}$.
By \refthm{ezremark} 
 there is symmetry between
$u_1v_1$ and $u_5v_5$; hence there are only two cases to consider,
namely $G_3=G_2+(u_1,v_1,u_2,u_3)$ and $G_3=G_2+(u_1,v_1,u_2,v_2)$,
where $G_2=G_1+(u_5,v_5,v_1',w)$, $G_4=G_3+(k_2,v_5,u_1,k_3)$,
$w$ is a neighbor of $v_1'$ different from $v_5$ and $k_i,l_i$
are the new vertices of $G_i$ for $i=2,3,4$. Let us first dispose
of the former case. Let $G',\eta'$ be obtained from $\eta$ by
first rerouting $u_1v_1$ along $u_5k_2k_4v_5$, and then
routing the new edge along $l_3k_3v_1$.
Then by \refthm{2.6} $G',\eta'$ satisfy the conclusion of the lemma
and so we may assume
that the latter case holds. We claim that we may assume that $w=v_4$.
Otherwise let $G'=G+(u_1,u_4,v_1',w)$ and $\eta'$ be obtained
from $\eta$ by routing the new edge along $u_5k_2l_2$; 
then $G',\eta'$ again satisfy the conclusion of the lemma by
\refthm{2.6}. Thus we may assume that $w=v_4$. Now $G_4$ is
isomorphic to a type F expansion of $G$,
and so the conclusion of the lemma is satisfied. \qed

\newthm{4.4} Let $F$ be a graph of minimum degree at least two, let
$G$ be a \qc\ cubic graph,  let $C$ be
a quadrangle in $G$, and let $G_5$ be a long $5$-extension of $G$
based at $C$ such that $F$ is a subgraph of both $G $ and $G_5$
and $F$ is disjoint from $C$.
Then there exist a graph $G'$ and a \he\ $\eta':G'\emb G_5$ such that
$\eta'$ fixes $F$ and $G'$ is either a type G or type H expansion
of $G$ based at $C$
or a long $n$-extension of $G$ based at $C$ for some $n\in\{1,2,3,4\}$.

\proof 
Similarly as in the previous two proofs we may assume that there exists
a short 1-extension $G_1=G+(u_1,u_4,v_1,v'_1)$ of
$G$ based at $C$ such that $G_5$ is a type F expansion of $G_1$ based at the new
quadrangle $C_1$ of $G_1$.
Since $G_1$ is a short extension, one of $u_1,u_4$ is adjacent to
one of $v_1,v_1'$, and so we may assume that $u_1$ is adjacent
to $v_1$. Let the vertices of $C$ be $u_1,u_2,u_3,u_4$ in order,
and for $i=1,2,3,4$ let $v_i$ be the neighbor of $u_i$ not on $C$.
Let $v_1''\not\in\{u_1,v_1'\}$ be the third neighbor of $v_1$.
Let $u_5,v_5$ be the new vertices of $G_1$; thus $V(C_1)= 
\{v_1,u_1,u_5,v_5\}$. 
Since by \refthm{ezremark} there is symmetry between
$u_1v_1$ and $u_5v_5$ there are only three cases to consider,
namely whether $G_5$ is based on $u_1u_5$, $v_1v_5$ or $u_5v_5$.
If $G_5$ is based on $u_1u_5$ then $G_5$ is a type G expansion
of $G$, and so  $G_5$ and the identity \he\ satisfy
the conclusion of the lemma. 

Next we assume that $G_5$ is based on $v_1v_5$.
It follows that $v_1'$ and $v_1''$ have a common neighbor, say $w$,
and $\{v_1,v_1',w,v_1''\}$ is the vertex-set of a quadrangle in $G$.
Let $z$ be the neighbor of $v_1'$ in $G$ other than $v_1$ and $w$.
By a rerouting argument similar to ones used in previous proofs it is easy
to construct a \he\ $G+(u_1,u_4,v_1',z)\emb G_5$ that fixes $F$.
By~\refthm{2.6} the lemma holds, unless $z=v_4$.
Thus we may assume that $z=v_4$, and similarly that $v_2$ and $v_1''$
have a common neighbor. We deduce that $G_5$ is a type H expansion
of $G$, as desired. 

We may therefore assume that $G_5$ is based on $u_5v_5$. Then
$v_4$ and $v_1'$ are adjacent. Let $G'$ be obtained from $G_5$
by deleting the edge $u_2u_3$ and suppressing the resulting degree
two vertices.
Then $G'$ is isomorphic to a type F expansion of $G$, 
and it is easy to construct a \he\ $\eta':G'\emb G_5$ that fixes $F$.
Hence $G',\eta'$ satisfy the conclusion of the lemma. \qed

\newthm{4.5}Let $F$ be a graph of minimum degree at least two, let
 $G$ be a \qc\ cubic graph, let $C$ be
a quadrangle in $G$, let $n\ge1$ be  an integer, and
let $G_3$ be a long $n$-extension of $G$ based at $C$
such that $F$ is a subgraph of both $G $ and $G_3$
and $F$ is disjoint from $C$.
Then there exist a graph $G'$ and a \he\ $\eta':G'\emb G_3$ such that
$\eta'$ fixes $F$ and $G'$ is a long $n'$-extension of $G$ based 
at $C$ for some $n'\in\{1,2,3,4,5\}$.

\proof
Similarly as in the previous three proofs we may assume that there exists
a short 1-extension $G_1=G+(u_1,u_4,v_1,v'_1)$ of
$G$ based at $C$ such that $G_3$ is a type G or H expansion of $G_1$ based at the new
quadrangle $C_1$ of $G_1$. 
Thus one of $u_1,u_4$ is adjacent to
one of $v_1,v_1'$, and so we may assume that $u_1$ is adjacent
to $v_1$. Let the vertices of $C$ be $u_1,u_2,u_3,u_4$ in order,
and for $i=1,2,3,4$ let $v_i$ be the neighbor of $u_i$ not on $C$.
Let $v_1''\not\in\{u_1,v_1'\}$ be the third neighbor of $v_1$ in $G$.
Let $u_5,v_5$ be the new vertices of $G_1$.
Thus  $V(C_1)=\{v_1,u_1,u_5,v_5\}$. Then
$G_1$ is \qc\ by \refthm{1extquadconn}.

We first assume that $G_3$ is a type G expansion of $G_1$.
Let $G_2$ be a $1$--extension of $G_1$ such that $G_3$ is a type
F expansion of $G_2$, and let $k_2,l_2$ be the new vertices of $G_2$.
\junk{and let $k_1,l_1,k_2,l_2,k_3,l_3,k_4,l_4$ be the new vertices
of $G_3$.} 
 From the symmetry it suffices to consider three subcases.
We consider them separately in the next three paragraphs.

As a first subcase assume that $G_2=G_1+(u_1,v_1,u_2,u_3)$.
Let $G',\eta'$ be obtained from $\eta$ by
first rerouting $\eta(u_1v_1)$ along $\eta(u_5v_5)$, 
and then routing the new edge along $l_2k_2v_1$.
Then $G',\eta'$ satisfy the conclusion of the lemma.

In the second subcase $G_2=G_1+(u_1,v_1,u_2,v_2)$, and $G_3$ is based
on $(u_1,k_2)$. Let $G'$  be obtained from $G_3$ by deleting the edge
$u_3u_4$, and suppressing the resulting vertices of degree two.
Then $G'$ is isomorphic to a type G expansion of $G$,
and it is easy to construct a  \he\ $\eta':G'\emb G_3$ that fixes $F$.
Then the pair $G',\eta'$ satisfies the conclusion of the lemma.

In the third and last subcase $G_2=G_1+(u_5,v_5,v_1',z)$, where
$z$ is a neighbor of $v_1'$ different from $v_5$, and $G_3$ is
based on $(v_5,k_2)$. Let $G'=G+(u_1,u_4,v_1',z)$;
by considering the path $u_5k_2l_2$ it is easy to construct a
\he\ $\eta':G'\emb G_3$  that fixes $F$.
If $z\not=v_4$ then by \refthm{2.6} the pair $G',\eta'$
satisfies the conclusion of the lemma. We may therefore assume
that $z=v_4$.
Let $G'$ be obtained from $G_3$ by deleting the edges $u_4v_4$ and
$u_2u_3$ and suppressing the resulting degree two vertices. 
and let
$\eta'$ be the canonical \he\ $G'\emb G_3$.
Then $G'$ is isomorphic to a type F expansion of $G$,
and  it is easy to construct a \he\ $\eta':G'\emb G_3$  that fixes $F$.
Hence the pair $G',\eta'$ satisfies the conclusion of the
lemma. 

We now assume that $G_3$ is a type H expansion of $G_1$ based at $C_1$.
By \refthm{ezremark} there is symmetry between $u_1$ and $u_5$.
Let $D$ be as in the definition of expansion of type H.
Thus some vertex of $D$ is adjacent in $G_1$ to some vertex of $C$.
By symmetry it suffices
to consider only two subcases.  
In the first subcase $v''_1\in V(D)$ is adjacent to $v_1\in V(C_1)$, and a
neighbor of $v''_1$ in $D$ is adjacent to $v'_1$.  But then the set
$V(D)\cup\{v_1,v'_1\}$ contradicts the \qcity\ of $G$.
In the second subcase some vertex of $D$ is adjacent to $u_5$, and 
the set $V(D)\cup\{u_1,u_2\}$ contradicts the \qcity\ of $G$.~\qed

\newthm{4.55}Let $F$ be a graph of minimum degree at least two, let
$G$ be a \qc\ cubic graph, let $C$ be
a quadrangle in $G$, let $n\ge1$ be an integer, and let $G_3$ be a long $n$-extension 
of $G$ based at $C$ such that $F$ is a subgraph of both $G $ and $G_3$
and $F$ is disjoint from $C$.
Then there exist a graph $G'$ and a \he\ $\eta':G'\emb G_3$ such that
$\eta'$ fixes $F$ and $G'$ is  a type A, B, C, D, E, F, G or H expansion
of $G$ based at $C$.

\proof
Let us choose an integer $n_2\ge1$, a graph $G_2$ and a \he\ $\eta_2:G_2\emb G_3$
such that $G_2$ is a long $n_2$-extension of $G$ based at $C$, the \he\  $\eta_2$ fixes $F$,
and, subject to that, $n_2$ is minimum.
Such a choice is possible, because $n_2=n$, $G_2=G_3$ and the identity \he\ 
satisfy the requirements (except minimality).

We claim that there do not exist an integer $n_1$, graph $G_1$ and
\he\ $\eta_1:G_1\emb G_2$ such that $1\le n_1<n_2$, $G_1$ is
a long $n_1$-extension of $G$ based at $C$ and $\eta_1$ fixes $F$.
Indeed, otherwise the graph $G_1$ and \he\ $\eta_1\circ\eta_2$
violate the choice of $G_2,\eta_2$.
This proves our claim that  $n_1, G_1,\eta_1$ do not exist.

It follows from~\refthm{4.1}, \refthm{4.2}, \refthm{4.3}, \refthm{4.4}, and \refthm{4.5}
that there exist a graph $G'$ and a \he\  $\eta':G'\emb G_2$ such that
$\eta'$ fixes $F$ and $G'$ is a  type A, B, C, D, E, F, G or H expansion
of $G$ based at $C$.
Thus $G'$ and the \he\ $\eta'\circ\eta_2$
satisfy the conclusion of the lemma.~\qed

\newthm{4.6}Let $G,H$ be cubic graphs, 
 let $F$ be a graph of minimum degree at least two,
let $\eta:G\emb H$ fix $F$, let $C$ be
a quadrangle in $G$ disjoint from $F$, let $G$
be \qc, and let $H$ be \cfc. Then
there exist an expansion $G'$ of $G$  based at $C$ of
type A, B, C, D, E, F, G, or H
and a homeomorphic embedding $\eta':G'\emb H$  that fixes~$F$.

\proof  
By \refthm{extex} there exist an integer $n\ge1$, a long $n$-extension
$G_2$ of $G$ based at $C$, and a \he\ $\eta_2:G_2\emb H$ that fixes $F$.
By \refthm{4.55} there exist a type A, B, C, D, E, F, G, or H expansion $G_1$ of $G$ 
based at $C$ and a \he\ $\eta_1:G_1\emb G_2$ that fixes $F$.
Thus $G_1$ and $\eta_1\circ\eta_2$ are as desired.~\qed

When $F$ is the null graph we obtain the following corollary.

\newthm{4.6cor}Let $G,H$ be cubic graphs, 
let $G\emb H$, let $C$ be
a quadrangle in $G$, let $G$
be \qc, and let $H$ be \cfc. Then
there exist an expansion $G'$ of $G$  based at $C$ of
type A, B, C, D, E, F, G, or H
and a homeomorphic embedding $G'\emb H$.

\newsection{dodec} DODECAHEDRAL CONNECTION

In this section we introduce \dcity, a notion of connectivity that is stronger than \cfcity. The main result of this section, \refthm{5.5} below, says that if the  graph $H$ in~\refthm{4.6} is \dc, then the last three outcomes of~\refthm{4.6} can be eliminated.

A \dfn{guild} is a pair $(G,\pi)$, where 
$G$ is a graph with every vertex of degree $1$ or $3$, and 
$\pi$ is a cyclic ordering of the set of vertices of $G$ of degree
$1$. (We consider $(1,2,3,4,5)$ and $(3,2,1,5,4)$ to be the same
cyclic ordering.) 
This is closely related to the notion of a society, 
introduced in~[\cite{RobSeyGM9}].
If $(G,\pi)$ and $(G',\pi')$ are guilds
and $\eta:G\emb G'$ is a \he, we say that $\eta$ is a \dfn{\he\ of
$(G,\pi)$ into $(G',\pi')$} if $\eta$ maps $\pi$ onto $\pi'$.
(That is, if $\pi=(v_1,v_2,\ldots,v_n)$, then $\pi'$ is the cyclic 
ordering $(\eta(v_1),\eta(v_2),\ldots,\eta(v_n))$.)
If that is the case we write $\eta:(G,\pi)\emb(G',\pi')$.
If $\delta A$ is an edge-cut of a cubic graph $G$ of cardinality
$k$ such that $\delta A$ is a matching, and $v_1,v_2,\ldots,v_k$ 
are all the vertices of $V(G)-A$
incident with an edge of $\delta A$, then let $H$ be the graph
$G\restriction (A\cup\{v_1,v_2,\ldots,v_k\})$. 
We say that
$(H,(v_1,v_2,\ldots,v_k))$ is a \dfn{shore guild corresponding 
to $A$}. Thus if $k>2$ there are $(k-1)!/2$ shore guilds
corresponding to $A$.

Let $G$ be the Dodecahedron, and let $C$ be a circuit of $G$
of length five with vertices $u_1,u_2,\ldots,u_5$ in order. 
For $i=1,2,\ldots,5$ let $v_i$ be the neighbor of $u_i$ not
on $C$. Let
$G'$ be the graph obtained from $G$ by deleting the edges of $C$;
then $D=(G',(u_1,u_2,u_3,u_4,u_5))$ is a guild, called the
\dfn{Dodecahedron guild}. Let $G''=G'+(u,v,x,y)$ be a $1$-extension
of $G'$. We say that $D'=(G'',(u_1,u_2,u_3,u_4,u_5))$ is 
a \dfn{non-planar expansion of the Dodecahedron guild} if
$\{u,v\}\not=\{u_i,v_i\}$ for all $i=1,2,\ldots,5$, and neither
$u$ nor $v$ is equal or adjacent to $x$ or $y$.

Let $G$ be a \cfc\ cubic graph. We say that $G$ is \dfn{\dc} if for every
edge-cut $\delta A$ of cardinality five and
every shore guild $S$ corresponding to $A$, if $\eta:D\emb S$
is a \he\ of the Dodecahedron guild into $S$, then there
exist a non-planar expansion $D'$ of $D$
and a \he\ $\eta':D'\emb S$.

The following proposition from [\cite{RobSeyThoCubic}] is not needed in this
paper, but is stated for the reader's convenience as it sheds some
light on the seemingly mysterious definition of dodecahedral
connection. A guild $(G,\pi)$ is \dfn{planar} if $G$ can be drawn
in a closed disc $\Delta$ with the vertices of degree one drawn in the
boundary of $\Delta$ in the order given by $\pi$.

\newthm{5.1}A \cfc\ cubic graph $G$ is \dc\ if and only if for every edge-cut
$\delta A$ of cardinality $5$ with $|A|\ge7$ and $|V(G)-A|\ge7$, 
no shore guild corresponding to $A$ is planar.


We need the following lemma.

\newthm{5.2}Let $G,G_1,H$ be cubic graphs, let $G$ be \qc,
let $F$ be a graph of minimum degree at least two,
let $C$ be a quadrangle
in $G$, let $G_1$ be a type F expansion of $G$ 
with core $R$ based at $C$ such that $R$ is disjoint from $F$, let $H$ be \dc, 
and let $\eta_1:G_1\emb H$ be a \he\ that fixes $F$. 
Then there exist a $1$-extension $G_2=G_1+(u,v,x,y)$
of $G_1$ and a homeomorphic embedding $G_2\emb H$ that fixes $F$ and such that
$u,v\in R$, and either $G_2$ is a long $1$-extension of $G_1$ or $x,y\not\in R$.

\proof Let $G,G_1,H,C,R,\eta_1$ be as stated. Then $\delta R$
is an edge-cut of $G_1$ of cardinality five such that some shore 
guild corresponding
to $R$ is isomorphic to the Dodecahedron guild. Let
$\delta R=\{e_1,e_2,\ldots,e_5\}$.
If there exists an edge-cut
$\delta A$ of $H$ of cardinality five with $\eta_1(R)\subseteq A$ and
$\eta_1(V(G_1)-R)\subseteq V(H)-A$ then the conclusion follows from
the definition of dodecahedral connection. We may therefore assume
that no such edge-cut exists. Thus by \refthm{3.3} there exists an augmenting
sequence $\gamma=(Q_1,Q_2,\ldots,Q_n)$ with respect to $G_1,H,R$ and $\eta_1$.
By \refthm{3.4} we may assume 
(by replacing $\eta_1$ by a different embedding if necessary)
that  the conclusion of \refthm{3.4} holds. Let $G_2,\eta_2$ be obtained from $\eta_1$
by routing the new edge along $Q_1$; it follows that $G_2$ and 
$\eta_2$ satisfy the conclusion of the lemma. \qed


The following result will allow us to eliminate type F expansions when
the graph $H$ is \dc.

\newthm{5.3}Let $G,G_4,H$ be cubic graphs, 
let $C$ be a quadrangle in $G$,
let $F$ be a graph of minimum degree at least two,
let $G$ be \qc, let $G_4$ be a type F expansion of $G$ based at $C$ with core $R$,
and let $G_5=G_4+(u,v,x,y)$ be a 1-extension of $G_4$ such that $u,v\in R$,
and either $G_5$ is a long 1-extension of $G_4$ or $x,y\not\in R$.
Assume further that $F$ is a subgraph of both $G$ and $G_5$.
Then there exist an integer $n\in\{1,2,3\}$, a long $n$-extension $G'$ of $G$
based at $C$,
and a homeomorphic embedding $\eta':G'\emb G_5$ that fixes $F$.

\proof
Let $u_1,u_2,u_3,u_4$ be the vertices
of $C$ in order, for $i=1,2,3,4$ let $v_i$ be the neighbor of
$u_i$ not on $C$, and let $v_i',v_i''$ be the neighbors of $v_i$
other than $u_i$. Since $G$ has a type F expansion 
we may assume that $v_1'=v_2'$. Let $w\not\in\{v_1,v_2\}$ be the
third neighbor of $v_1'$.
Choose $G_1,G_2,G_3$ such that each of $G,G_1,G_2,G_3,G_4,G_5$
 is a $1$--extension of the previous. 
For $i=0,1,2,3,4$ let $\eta_i$ be the canonical \he\ $G_i\emb G_5$
determined by the generating sequence $G_{i+1},G_{i+2},\ldots,G_5$,
where $G_0$ means $G$,
for $i=1,2,3,4$ let 
let $k_i,l_i$ be the new vertices of $G_i$, and let
$G_1=G+(u_1,u_2,v_1,v_1')$, $G_2=G_1+(k_1,u_2,v_2,v_1')$,
$G_3=G_2+(k_1,k_2,u_3,u_4)$ and $G_4=G_3+(k_1,l_1,k_2,l_2)$.
Then $R=\{u_1,u_2,u_3,u_4,v_1,v_2,v_1',k_1,l_1,k_2,l_2,k_3,k_4,l_4\}$.
 From \refthm{2.6} we deduce that

\claim{1}the vertices in $R\cup\{v_1'',w,v_2'',v_3,v_4\}$ are pairwise
distinct, except that possibly $w=v_3$ or $w=v_4$, but not both.

We also point out for future reference that

\claim{2}there is symmetry fixing $v_1',w$ and taking 
$u_1,u_4,v_1,v_1'',v_4$ onto $u_2,u_3,v_2,v_2'',v_3$, respectively.

\claim{3}If $u\in\{k_1,k_2,k_3,k_4,l_4\}$ and $x,y\not\in R-\{v_1'\}$ 
then the lemma holds.

To prove (3) let $u,x,y$ be as stated, and let $G'=G+(u_1,u_2,x,y)$
and $\eta'$ be obtained from $\eta_0$ by routing the new edge along
$\eta_5(k_5l_5)\cup Q$, where $Q$ is an appropriate subpath of
$\eta_4(G_4)$. Then $G',\eta'$ satisfy the conclusion of the lemma,
and (3) follows.
\medskip

\claim{4}If $u=l_3$,  $\{x,y\}\cap\{v_3,v_4\}=\emptyset$, $\{x,y\}\not=
\{u_1,k_1\}$ and $\{x,y\}\not=\{u_2,k_2\}$ then the lemma holds.

To prove (4) we may assume by (2) that
$\{x,y\}$ does not equal $\{u_2,v_2\}$, $\{k_1,k_4\}$ or 
$\{k_4,l_1\}$. Let
$G',\eta'$ be obtained from $\eta_0$ by rerouting $\eta_5(u_1v_1)$
along $\eta_1(k_1l_1)$, and then 
by routing the new edge along
$Q\cup\eta_5(k_5l_5)\cup Q'$, where $Q$ is $\eta_5(k_5l_3)$ if
$v=k_3$ and null otherwise, 
and $Q'$ is $\eta_5(l_5v_1)$ if $\{x,y\}=\{u_1,v_1\}$, a subpath
of $\eta_4(k_2l_4)\cup\eta_4(l_4l_2)\cup\eta_4(l_4k_4)$ with
ends $\eta_5(l_5)$ and $\eta_4(l_2)$ if $l_4\in\{x,y\}$, and null
otherwise.
Then $G',\eta'$ satisfy the conclusion of the lemma and (4)
follows.
\medskip

\claim{5}If $\{u,v\}=\{k_3,l_3\}$, then the lemma holds.

This follows immediately from (3) and (4).
\medskip

\claim{6}If $u\in\{u_2,k_2\}$ and  $x,y\not\in\{k_1,v_2,l_2,u_3,
l_3,k_4,l_4,v_2''\}$ then the lemma holds.

To prove (6) let $G'$ and $\eta'$ be obtained from $\eta_0$ by
first rerouting $\eta_0(u_2u_3)$
along $\eta_5(k_3l_3)$, then routing the first new edge along
$\eta_5(k_1k_4)\cup\eta_5(k_4l_4)\cup\eta_5(l_4l_2)$, and then
routing the second new edge along $\eta_5(k_5l_5)\cup Q$, where
$Q$ is $\eta_5(k_5u_2)$ if $v=u_3$, $\eta_5(k_5k_2)$ if $v=l_4$,
and null otherwise. Then $G',\eta'$ satisfy
the conclusion of the lemma, because $G'$ is a long $2$-extension of $G$.
This proves (6).
\medskip

\claim{7}If $v_3$ and $v_2''$ are adjacent, then the lemma holds.

To prove (7) let $G',\eta'$ be 
obtained from $\eta_0$ by 
\item{$\bullet$}
first rerouting $\eta_4(u_2u_3)$ along $\eta_4(k_3l_3)$, 
\item{$\bullet$}
then rerouting $\eta_4(k_2u_2)\cup\eta_4(u_2v_2)$
along $\eta_2(k_2l_2)$, 
\item{$\bullet$}
then rerouting $\eta_4(u_1v_1)$ along $\eta_1(k_1l_1)$, 
\item{$\bullet$}
then rerouting $\eta_4(k_1k_3)$ along $\eta_4(k_4l_4)$, 
\item{$\bullet$}
then rerouting $\eta_4(v_3v_2'')$ along $\eta_4(u_3u_2)\cup\eta_4(u_2v_2)$, 
\item{$\bullet$}
then routing the first new edge along $\eta_4(u_2k_2)$, 
\item{$\bullet$}
then routing the second new edge along $\eta_4(k_1k_3)$, 
\item{$\bullet$}
and finally routing the third new edge along $\eta_4(u_1v_1)$.

\noindent
Then $G'$ is a type E expansion of $G$, 
and hence the pair $G',\eta'$ satisfies
the conclusion of the lemma. This proves (7).
\medskip

\claim{8}If  $\{u,v\}$ is one of $\{u_2,u_3\}$,
$\{u_3,l_3\}$, $\{l_3,u_4\}$ or 
$\{u_1,u_4\}$, then the lemma holds.

To prove (8) we may assume by (2) that $\{u,v\}=\{u_2,u_3\}$ or
$\{u_3,l_3\}$. Assume first that $v_3,v_4\not\in\{x,y\}$. Let
$G',\eta'$ be obtained from $\eta_0$ by 
\item{$\bullet$}
first rerouting
$\eta_4(u_1u_4)$ along $\eta_4(l_3k_3)$, 
\item{$\bullet$}
then rerouting
$\eta_4(k_3k_1)\cup\eta_4(k_1u_1)\cup\eta_4(u_1v_1)$ along $\eta_4(k_2l_4)\cup
\eta_4(l_4k_4)\cup\eta_4(k_4l_1)$, 
\item{$\bullet$}
then rerouting
$\eta_4(l_1v_1')$ along $\eta_4(l_4l_2)$ and 
\item{$\bullet$}
finally
routing the new edge
along $\eta_5(k_5l_5)\cup Q$, where $Q$ is either null, or a  path of
$\eta_4(G_4)$ with one end $\eta_5(l_5)$, the other end in
$\eta'(v_1v_1'')$, and otherwise disjoint from $\eta'(G)$. 

\noindent
The graph $G'$ is a long extension of $G$, unless $\{u,v\}=\{u_3,l_3\}$ and
$\{x,y\}=\{k_2,l_4\}$, in which case (8) follows from (4). 
Thus (8) holds if $v_3,v_4\not\in\{x,y\}$, and so we may assume that
either $x=v_3$ or $x=v_4$. As a second case assume that $x=v_4$. 
If $\{u,v\}=\{u_2,u_3\}$, then
(8) follows from (6), and so let $\{u,v\}=\{u_3,l_3\}$.
Let $G'$ be obtained from $G_4$ by deleting the edges $k_2l_4$ and
$k_4l_1$ and suppressing degree two vertices. Then $G'$ is
isomorphic to a type E expansion of $G$, and so (8) follows. This
completes the second case. Thirdly, let $x=v_3$. Since the cases
$\{u,v\}=\{u_2,u_3\}$ and $\{u,v\}=\{u_3,l_3\}$ are symmetric by \refthm{ezremark},
we may assume that $\{u,v\}=\{u_2,u_3\}$. If $v_3$ and $v_2''$
are adjacent, then (8) follows from (7); otherwise it follows
from (6). This proves (8).
\medskip

\claim{9}If $u=k_3$ then the lemma holds.

To prove (9) let $u=k_3$. By (5) we may assume that $v\not=l_3$
(and hence $\{x,y\}\ne\{k_4,l_4\}$),
by (3) we may assume that $\{x,y\}\cap R\not=\emptyset$,
and by (2) we may assume that $\{x,y\}\not=\{u_2,v_2\}$ and
$\{x,y\}\not=\{l_1,k_4\}$. By (8) we may assume that 
$\{x,y\}\not=\{u_2,u_3\}$ and $\{x,y\}\not=\{u_1,u_4\}$.
Let $G',\eta'$ be obtained from $\eta_0$ first by rerouting
$\eta_4(k_1k_3)\cup\eta_4(k_3k_2)$ along $\eta_4(k_1k_4)\cup
\eta_4(k_4l_4)\cup\eta_4(l_4k_2)$, then rerouting $\eta_4(u_1v_1)$
along $\eta_4(k_4l_1)$, and finally routing the new edge along
$\eta_4(l_3k_3)\cup \eta_5(k_3k_5)\cup\eta_5(k_5l_5)\cup Q$, where $Q$
is either null or $\eta_5(l_2l_5)$ or $\eta_5(v_1l_5)$.
If $\{x,y\}\not
=\{u_4,v_4\}$ and $\{x,y\}\not=\{u_3,v_3\}$ then $G'$ is a long extension of $G$,
and hence the lemma holds. From the symmetry we may assume that
$\{x,y\}=\{u_4,v_4\}$. If $\{u,v\}=\{k_3,k_2\}$ then (9) follows
from (6), and so we may assume that $\{u,v\}=\{k_1,k_3\}$.
Let $G',\eta'$ be obtained from $\eta_0$ by first rerouting
$\eta_4(u_1u_4)$ along $\eta_5(k_5l_5)$, then rerouting 
$\eta_4(k_1u_1)\cup\eta_4(u_1v_1)$ along $\eta_4(k_1k_4)\cup
\eta_4(k_4l_1)$, and then routing the new edge along $\eta_4(u_4u_1)
\cup\eta_4(u_1v_1)$. Then $G',\eta'$ satisfy the conclusion of the
lemma, and hence (9) holds.
\medskip

\claim{10}If $u\in\{l_1,k_4\}$ and $x=v_2$ then the lemma holds.

To prove (10) we first define two paths $Q,Q'$. Let $Q$ be the
path of $\eta_4(v_2l_2)\cup\eta_4(v_2u_2)\cup\eta_4(v_2v_2'')$
with one end $\eta_5(l_5)$ and the other end in
$\eta_4(v_2u_2)\cup\eta_4(v_2v_2'')$, and let $Q$ be the path of 
$\eta_1(k_1l_1)\cup\eta_0(v_1v_1')\cup\eta_4(k_4l_4)$ with one
end $\eta_5(k_5)$ and the other end in $\eta_0(v_1v_1')$. Let
$G'=G+(u_1,u_2,v_1',w)$ and $\eta'$ be obtained from $\eta_0$
by rerouting an appropriate path along $Q\cup \eta_5(l_5k_5)\cup Q'$,
and then routing the new edge along $\eta_2(k_2l_2)\cup 
\eta_2(l_2v_1')$. Then $G',\eta'$ satisfy the conclusion of the 
lemma, thus proving (10).
\medskip

\claim{11}If $\{u,v\}$ equals one of $\{u_2,v_2\}$, $\{u_1,v_1\}$,
$\{u_2,k_2\}$ or $\{u_1,k_1\}$, then the lemma holds.

To prove (11) we may assume by (2) that $\{u,v\}=\{u_2,v_2\}$
or $\{u,v\}=\{u_2,k_2\}$. If $\{x,y\}=\{u_4,v_4\}$ or
$\{x,y\}=\{v_3,v_3'\}$, where $v_3'\not=u_3$ is a neighbor of $v_3$,
then (11) follows from (6) and (7). If $\{x,y\}=\{k_1,k_3\}$ then (11)
follows from (9), and if $\{x,y\}=\{k_4,l_4\}$, then (11) follows
from (10). We may therefore assume that none of the above hold.
Let $G',\eta'$ be obtained from $\eta_0$ by first rerouting
$\eta_4(u_1u_4)$ along $\eta_4(k_3l_3)$, then rerouting
$\eta_4(v_1l_1)$ along $\eta_1(k_1l_1)$, then rerouting
$\eta_4(u_2v_2)$ along $\eta_2(k_2l_2)$, then rerouting
$\eta_4(l_2v_1')$ along $\eta_4(k_4l_4)$, and finally routing
the new edge along $Q\cup\eta_5(k_5l_5)\cup Q'$, where $Q$ is either
null or $\eta_5(k_5u_2)$, and $Q'$ is either null or a subpath
of $\eta_4(l_2v_1')\cup\eta_4(u_4u_1)\cup\eta_4(v_1l_1)$ with
one end $\eta_5(l_5)$ and the other end in $\{\eta_4(v_1'),\eta_4(u_1),
\eta_4(v_1)\}$. Then the graph $G'$ is a long extension of $G$, 
and hence (11) holds.
\medskip

\claim{12}If $u\in\{l_1,k_4\}$, $x\in V(G)-(R\cup\{v_1'',w\})$ and
$y\not=v_1''$, then the lemma holds.

To prove (12) let $G',\eta'$ be obtained from $\eta_0$ first by
rerouting $\eta_4(u_4u_1)$ along $\eta_4(k_3l_3)$, then rerouting
$\eta_4(k_1k_3)$ along $\eta_4(k_2l_4)\cup\eta_4(l_4k_4)\cup
\eta_4(k_4k_1)$, then by rerouting $\eta_0(v_1v_1')$ along
$\eta_4(l_4l_2)$, then routing a first new edge along 
$\eta_4(k_1k_3)$, and finally  routing a second new edge along
$Q\cup\eta_5(k_5l_5)$, where $Q$ is a suitable path of $\eta_4(G_4)$
with one end $\eta_5(k_5)$ and the other end in $\eta_4(k_1k_4)\cup
\eta_4(k_4l_4)$. Then $G',\eta'$ satisfy the conclusion of (12)
thus proving (12).
\medskip

\claim{13}If $\{u,v\}=\{k_4,l_4\}$ then the lemma holds.

This follows from (2), (3), (5), (8), (10), (11) and (12).
\medskip

\claim{14}If $\{u,v\}$ equals one of $\{k_1,k_4\}$, $\{k_4,l_1\}$,
$\{k_2,l_4\}$ or $\{l_4,l_2\}$ then the lemma holds.

To prove (14) we may assume by (2) that $\{u,v\}=\{k_1,k_4\}$ or
$\{u,v\}=\{k_4,l_1\}$.
By (3), (5), (8), (10), (11) and (12) we may assume that $\{u,v\}=
\{k_1,k_4\}$, and
$\{x,y\}=\{l_2,v_1'\}$ or $\{x,y\}=\{v_1,v_1''\}$. Let $\{u,v\}=
\{k_1,k_4\}$, and assume first that
$\{x,y\}=\{l_2,v_1'\}$. Let $G',\eta'$ be obtained from $\eta_0$
first by rerouting $\eta_0(v_1'l_2)$ along $\eta_4(l_1k_4)\cup
\eta_4(k_4l_4)\cup\eta_4(l_4l_2)$, and then routing the new edge along
$\eta_5(k_1k_5)\cup\eta_5(k_5l_5)\cup\eta_5(l_5v_1')$. Then $G',\eta'$
satisfy the conclusion of the lemma. We may therefore assume that
$\{x,y\}=\{v_1,v_1''\}$. In this case let $G',\eta'$ be obtained
from $\eta_0$ by routing the first new edge along
$\eta_4(k_1k_4)\cup\eta_4(k_4l_4)\cup\eta_4(l_4l_2)$, and routing
the second new edge along $\eta_5(k_5l_5)$. Then $G'$ is a long
$2$-extension of $G$ by
\refthm{1extquadconn} and \refthm{2.4} (or by (1)), and hence the pair $G',\eta'$ satisfies
the conclusion of the lemma. This proves (14).
\medskip

\claim{15}If $\{u,v\}=\{v_2,l_2\}$ or $\{u,v\}=\{v_1,l_1\}$ then
the lemma holds.

To prove (15) we may assume by (2) that $\{u,v\}=\{v_1,l_1\}$.
By (5), (8), (9), (10), (11) and (12) we may assume that $\{x,y\}=
\{v_1'',z\}$, where $z\not=v_1$ is a neighbor of $v_1''$. But
then $G_5$ is isomorphic to $G_4+(u_1,v_1,v_1'',z')$, where $z'\not
\in\{v_1,z\}$ is the third neighbor of $v_1''$, and hence (15)
follows from (11).
\medskip

\claim{16}If $\{u,v\}=\{l_2,v_1'\}$ or $\{u,v\}=\{l_1,v_1'\}$ then
the lemma holds.

To prove (16) we may assume by (2) that $\{u,v\}=\{l_2,v_1'\}$.
By (2), (5), (8), (9), (11), (12) and (14) we may assume that 
$\{x,y\}=\{v_2'',z\}$, where $z\not=v_2$ is a neighbor 
of $v_2''$. 
Let $G',\eta'$
be obtained from $\eta_0$ by first rerouting $\eta_4(v_2v_2'')$
along $\eta_5(k_5l_5)$, then rerouting $\eta_4(u_2v_2)\cup
\eta_4(v_2l_2)$ along $\eta_2(k_2l_2)$, and finally routing the 
new edge along $\eta_4(u_2v_2)\cup\eta_4(v_2v_2'')$. If $v_2''$
is not adjacent to $v_3$, then $G'$ is a long extension of $G$ by 
\refthm{1extquadconn}, and hence
the pair $G',\eta'$ satisfies the conclusion of the lemma. On the
other hand if $v_3$ and $v_2''$ are adjacent, then (16) follows
from (7). This completes the proof of (16).
\medskip

The lemma now follows from (5), (8), (9), (11), (13), (14),
(15) and (16). \qed


\newthm{5.35}Let $G,G_4,H$ be cubic graphs, let $C$ be a quadrangle in $G$,
let $F$ be a graph of minimum degree at least two,
let $G$ be \qc, let $G_4$ be a type F expansion of $G$ based at $C$
such that its core is disjoint from $F$,
let $\eta:G_4\emb H$ fix $F$, and let $H$ be \dc. Then
there exist an integer $n\in\{1,2,3\}$,
a long $n$-extension $G'$ of $G$ based at $C$
and a homeomorphic embedding $\eta':G'\emb H$ that fixes $F$.

\proof
This follows immediately from~\refthm{5.2} and~\refthm{5.3}.~\qed

\newthm{5.4}Let $G,G_5,H$ be cubic graphs, let $C$ be a quadrangle in $G$,
let $F$ be a graph of minimum degree at least two,
let $G$ be \qc, let $G_5$ be a type G or H expansion of $G$ based at $C$,
let $F$ be a subgraph of both $G$ and $G_5$,
and let $H$ be \dc. Then
there exist an integer $n\in\{1,2,3\}$,
a long $n$-extension $G'$ of $G$ based at $C$
and a homeomorphic embedding $\eta':G'\emb H$ that fixes $F$.

\proof 
Let $G_1$ be a short $1$--extension of $G$ such that $G_5$ is a
type F expansion of $G_1$ based at the new quadrangle of $G_1$. 
By \refthm{5.35} applied to $G_1$ and
the new quadrangle of $G_1$ there exist an integer $k\in\{1,2,3\}$,
a long $k$-extension $G_2$ of $G_1$, and a \he\ $G_2\emb H$.
Then $G_2$ is a long $(k+1)$-extension of $G$ based at $C$,
and  so if $k\le2$, then the lemma holds.
We may therefore assume that $k=3$.
By~\refthm{4.3} we may assume that there exist a
type F expansion $G_3$ of $G$ based at $C$ and a \he\ $G_3\emb H$
that fixes $F$.
The conclusion of the lemma now follows from~\refthm{5.35}
applied to the graph $G$ and quadrangle $C$.~\qed

\newthm{5.5}Let $G,H$ be cubic graphs, 
let $C$ be a quadrangle in $G$, 
let $F$ be a graph of minimum degree at least two disjoint from $C$,
let $\eta:G\emb H$ fix $F$, 
let $G$ be \qc, and let $H$ be \dc. 
Then there exist an expansion $G'$ of $G$ of type A, B, C, D, or E based
at $C$,
and a homeomorphic embedding $G'\emb H$ that fixes $F$.

\proof By \refthm{4.6} there exist an expansion $G_1$ of $G$ of type
A, B, C, D, E, F, G or H and a \he\ $\eta_1:G_1\emb H$ that fixes $F$.
We may assume that $G_1$ is of type F, G, or H,
for otherwise $G_1,\eta_1$ satisfy the theorem.
By \refthm{5.35} and \refthm{5.4} applied to $G,G_1,H$ and $\eta_1$ 
there exist an integer $n\in\{1,2,3\}$, a long $n$-extension
$G_2$ of $G$ based at $C$ and a \he\ 
$\eta_2:G_2\emb H$  that fixes $F$. 
By \refthm{4.1} and \refthm{4.2}
 there exist an expansion $G_3$ of $G$ of type A, B, C, D, or E
and a \he\ $\eta_3:G_3\emb H$ that fixes $F$, as desired.~\qed

\newsection{two} A TWO-EXTENSION THEOREM

In this section we prove a preliminary  weaker version of~\refthm{1.3}.
In~\refthm{6.1} we prove it when $H$ is \dc, and in~\refthm{6.2}
we prove it for \cfc\ graphs~$H$.

\newthm{6.1}Let $G,H$ be cubic graphs, let $G$ be \cfc, let $H$
be \dc, and let $\eta:G\emb H$
be a homeomorphic embedding. Then there exist a \cfc\ cubic graph 
$G'$ and a homeomorphic embedding $\eta':G'\emb H$ such that 
$G'$ is  a $1$- or $2$-extension of $G$.

\proof Let $G,H,\eta$ be as stated. By \refthm{3.6} there exist
a $1$--extension $G_0=G+(u_2,v_3,u_1,v_4)$ of $G$ and a 
homeomorphic embedding $\eta_0:G_0\emb H$.
Let $u_3,u_4$ be the new vertices of $G_0$.
If $G_0$ is \cfc, then $G_0,\eta_0$ satisfy
the conclusion of \refthm{6.1}, and so we may assume that $G_0$ is
not \cfc. By \refthm{1extquadconn} we may assume that say $u_1$ is adjacent to
$u_2$. Then $G_0$ is \qc, and has a unique quadrangle $C_0$,
where $V(C_0)=\{u_1,u_2,u_3,u_4\}$. By \refthm{5.5} there exist an
expansion $G_2$ of $G_0$  of type A, B, C, D or
E based $C_0$, and a homeomorphic embedding $\eta_2:G_2\emb H$.
If $G_2$ is of type A, then the pair $G_2,\eta_2$ satisfies the conclusion of the lemma,
and so it remains to consider types B, C, D and E.

Let us assume now that $G_2$ is of type B, C or D, and let $G_1,G_2$ 
be a standard generating sequence for $G_2$.
Let $G_1=G_0+(a_1,a_2,v_1,w)$, where  $a_1,a_2\in V(C_0)$, 
$v_1,w\not\in V(C_0)$, and let $k_1,l_1$ be the new vertices of $G_1$.
%
Since $G_1$ is a short extension of $G_0$,
by \refthm{1extquadconn} we may assume that say $a_1$ is
adjacent to $v_1$, and hence $G_1$ has a unique quadrangle, say $C_1$,
and its vertex-set is $\{a_1,v_1,k_1,l_1\}$.
Let $\xi:G\emb G_2$ be the canonical \he\ determined by the generating
sequence $G_0,G_1,G_2$,  let
$\zeta_0=\xi\circ\eta_2$,
 let $\zeta_1$
be obtained from $\zeta_0$ by rerouting $\zeta_0(u_1u_2)$ along
$\eta_2(u_3u_4)$, and let $\zeta_2$ be obtained from $\zeta_0$
by rerouting $\eta_2(a_1v_1)$ along $\eta_2(k_1l_1)$.

\claim{1}We may assume that $\{a_1,a_2\}=\{u_1,u_2\}$.

To prove (1) we first notice that by \refthm{ezremark} we may assume
that $\{a_1,a_2\}=\{u_1,u_2\}$ or $\{a_1,a_2\}=\{u_3,u_4\}$. 
But if $\{a_1,a_2\}=\{u_3,u_4\}$, then by
replacing $\eta$ by $\zeta_1$ 
we can arrange that (1) holds. 
%
\medskip

 From the symmetry between $u_1$ and $u_2$ we may assume that
$a_1=u_1$ and $a_2=u_2$. Then $v_1$ is the neighbor of $u_1$ in $G_0$
that does not belong to $C_0$. Let $v_2$ be the neighbor
of $u_2$ in $G_0$ that does not belong to $C_0$.

\claim{2}We may assume that $w$ and $v_2$ are adjacent in $G$.

To prove (2) suppose that $w$ and $v_2$ are not adjacent,
and let $G'=G+(u_2,v_2,v_1,w)$ and $\eta'$ be obtained from $\zeta_1$
by routing the new edge along
$\zeta_1(u_2k_1)\cup\zeta_1(k_1l_1)$. Since 
$G+(u_2,v_2,v_1,w)$ is \cfc\ by  \refthm{2.2},  (2) holds.
\medskip

By (2) $G$ has a circuit with vertex-set $\{u_1,u_2,v_2,w,v_1\}$. 
Let $v_1'$ be the
neighbor of $v_1$ not on this circuit, and let $v_2'$ be defined
similarly. We distinguish cases depending on the type of the
expansion $G_2$.

Let us assume first that $G_2$ is a type B expansion of $G_0$.
Then $G_2=G_1+(v_1,l_1,x,y)$ for some $x,y\in V(G_1)$. Let
$k_2,l_2$ be the new vertices of $G_2$.
Let us assume first that $\{x,y\}=\{u_4,v_4\}$. Let 
$G'=G+(u_2,v_2,v_1,v_1')$ and $\eta'$ be obtained from $\zeta_1$
by first rerouting $\zeta_1(v_1u_1)$ along $\eta_2(k_2l_2)$,
and then routing the new edge along $\eta_1(u_1u_2)\cup\eta_1(u_1v_1)$.
Since $G+(u_2,v_2,v_1,v_1')$ is \cfc\ by \refthm{2.4}, the pair
$G',\eta'$ satisfies the conclusion of the theorem, as required.
We may therefore assume that $\{x,y\}\not=\{u_4,v_4\}$. 
Let $G',\eta'$ be obtained from $\zeta_0$ by routing the first
new edge along $\eta_1(k_1l_1)$, and then routing the second new edge
along $\eta_2(k_2l_2)$ (or along $\eta_2(k_2l_2)\cup\eta_2(l_2u_3)$
if $\{x,y\}=\{u_3,u_4\}$). Then $G',\eta'$ satisfy the conclusion
of the theorem.
This completes the case when $G_2$ is a type B
expansion of $G_0$.

We now assume that $G_2$ is a type C expansion of $G_0$. Since
$G_1+(k_1,l_1,v_2,w)$ is not \cfc, there are only two cases to
consider. Assume first that $G_2=G_1+(k_1,l_1,v_2,v_2')$,
and let $k_2,l_2$ be the new vertices of $G_2$.
Let $G'=G+(u_1,v_1,v_2,v_2')$ and $\eta'$ be obtained from $\zeta_2$
by routing the new edge along $\eta_2(k_2l_2)$.
Since $G'$ is \cfc\ by \refthm{2.4} the
theorem holds. Secondly, let us assume that 
$G_2=G_1+(u_1,v_1,v_4,v_4')$, where $v_4'\not=u_4$ is a neighbor 
of $v_4$ in $G$, and let $k_2,l_2$ be the new vertices of $G_2$.
Let $G'=G+(v_1,v_1',v_4,v_4')$ and $\eta'$ be obtained from
$\zeta_2$ by routing the new edge along 
$\eta_2(v_1k_2)\cup\eta_2(k_2l_2)$.
If $G'$ is \cfc, then the pair $G',\eta'$ is as desired.
We may therefore assume that 
$G+(v_1,v_1',v_4,v_4')$ is not \cfc, and hence $v_1'$ and $v_4'$
are adjacent by \refthm{2.2}. 
Let $G'$ and $\eta'$ be obtained from 
$\zeta_2$ by first rerouting $\zeta_2(v_1'v_4')$ 
along $\eta_2(v_1k_2)\cup\eta_2(k_2l_2)$, then  routing a first new
edge along $\eta_1(u_3u_4)$ and then routing a second new edge
along $\eta_2(k_2u_1)$. Then $G'$ is isomorphic to
$G_0+(u_1,u_4,v_1',v_4')$. Since
$G_0+(u_1,u_4,v_1',v_4')$ is \cfc\ by \refthm{2.2}, the pair $G',\eta'$
is as desired.
This completes the case when $G_2$ is a type C expansion.

We now assume that $G_2$ is a type D expansion of $G_0$; then
$G_2=G_1+(k_1,u_1,u_3,v_3)$. Let $k_2,l_2$ be the new vertices of 
$G_2$. Let $G'=G+(v_1,w,u_2,v_3)$ and $\eta'$ be obtained from
$\zeta_1$ by routing the new edge along
$\eta_2(k_2l_2)\cup\eta_2(k_2k_1)\cup\eta_2(k_1l_1)$.
Since $G'$ is \cfc\ by
\refthm{2.4} the theorem holds in this case. This completes the case
that $G_2$ is a type D expansion.

Finally we assume that $G_2$ is a type E expansion of $G_0$. Let
$G_1,G_2',G_2$ be a standard generating sequence for $G_2$.
From the symmetry we may assume that
$G_2'=G_1+(u_3,u_4,v_3,v_3')$, where $v_3'\not=u_3$ is a neighbor
of $v_3$,
and $G_2=G_2'+(k_2',u_3,k_1,u_1)$, where
$k_2',l_2'$ are the new vertices of $G_2'$. 
 Let $k_2,l_2$ be the new vertices of 
$G_2$. Let $G'$ and $\eta'$ be obtained from $\zeta_0$ by routing
the first new edge along 
$\eta_2(l_2k_2)\cup\eta_2(k_2k_2')\cup\eta_2(k_2'l_2')$ 
and then routing the second new edge along $\eta_2(k_1l_1)$.
Since $G'$ is \cfc\ by \refthm{2.3},
the theorem holds in this case. This completes the case when $G_2$
is a type E expansion of $G_0$, and hence the proof of the theorem. 
\qed

Let us recall that circuit expansion was defined prior to (1.3).

\newthm{6.2}Let $G,H$ be non-isomorphic \cfc\ cubic graphs, and 
let $\eta:G\emb H$
be a homeomorphic embedding. Then there exist a \cfc\  cubic graph 
$G'$ and a homeomorphic embedding $\eta':G'\emb H$ such that 
$G'$ is either a $1$- or $2$-extension or a 
circuit expansion of $G$. 

\proof Let $G,H,\eta$ be as stated. By \refthm{3.6} there exist
a $1$--extension $G_0=G+(u_2,v_3,u_1,v_4)$ of $G$ and a 
homeomorphic embedding $\eta_0:G_0\emb H$.
Let $u_3,u_4$ be the new vertices of $G_0$.
If $G_0$ is \cfc, then $G_0,\eta_0$ satisfy
the conclusion of \refthm{6.1}, and so we may assume that $G_0$ is
not \cfc. By \refthm{1extquadconn} we may assume that say $u_1$ is adjacent to
$u_2$. Then $G_0$ is \qc, and has a unique quadrangle $C_0$,
where $V(C_0)=\{u_1,u_2,u_3,u_4\}$. By \refthm{4.6cor} there exist an
expansion $G_2$ of $G_0$  of type A, B, C, D,
E, F, G or H based at $C_0$, and a homeomorphic embedding $\eta_2:G_2\emb H$.
If  $G_2$ is an expansion
of type A, B, C, D or E then the theorem holds by the proof of
\refthm{6.1}. We may therefore assume that $G_2$ is an expansion of
type F, G or H. Let $\zeta_1$ be defined as in the proof of \refthm{6.1}.

Assume first that $G_2$ is an expansion of type F. Since $G$ is
\cfc, $v_1$ and $v_4$ have no common neighbor in $G$, and
similarly $v_2$ and $v_3$ have no common neighbor in $G$.
Therefore $G_2$ is based on either $u_1u_2$, or $u_3u_4$.
In either case $G_2$ is a circuit expansion of $G$, and so the pair $G_2,\eta_2$
satisfies the conclusion of the theorem.

Secondly, let us assume that $G_2$ is an expansion of type G.
Let $G_1$ be a short $1$--extension of $G_0$ based at $C_0$
such that $G_2$ is a type F expansion of $G_1$, and let $C_1$ be the
unique quadrangle of $G_1$.
By replacing $\eta$ by $\zeta_1$ and by using symmetry
we may assume that $G_1=G_0+(u_1,u_2,v_1,w)$, where $w\not=u_1$
is a neighbor of $v_1$. Let $k_1,l_1$ be the new vertices of $G_1$;
then the vertex-set of $C_1$ is $\{u_1,v_1,l_1,k_1\}$. From claim
(1) in the proof of \refthm{6.1} we may assume that $w$ and $v_2$ are 
adjacent. Let $G'$ be obtained from $G_2$ by deleting the edge
$v_2w$ and suppressing the resulting vertices of degree two,
and let $\eta'$ be the restriction of $\eta_2$ to $G'$. Then $G'$
is isomorphic to a circuit expansion of $G$, and so the theorem
holds. 

Finally let us assume that $G_2$ is an expansion of type H.
Using the same symmetry as before we may assume that
$G_0$ has a quadrangle $D$ with vertex-set $\{x_1,x_2,x_3,x_4\}$,
where $u_1$ is adjacent to $x_1$, the vertices $u_2$ and $x_2$ have a
common neighbor, and $u_4$ and $x_4$ have a
common neighbor, say $z$.
Then the set $V(D)\cup\{u_1,z\}$ violates the \dcity\ of $G$.
This completes the case when $G_2$  is a type H expansion, 
and hence a proof of the theorem.~\qed

\newsection{one} A ONE-EXTENSION THEOREM

In this section we prove \refthm{1.3} and \refthm{1.4}.

\newthm{7.1}Let $G,H$ be \cfc\ cubic graphs, let $u_1,u_2,u_3,u_4,u_5$
(in order) be the vertices of a path of $G$, 
let $G_2=G+(u_1,u_2,u_3,u_4)+(u_2,u_3,u_4,u_5)$, 
and let
$\eta_2:G_2\emb H$.
Then there exist a \cfc\ handle expansion $G'$ of $G$ and
a homeomorphic embedding $G'\emb H$.
 
\proof Let $v_2\not\in\{u_1,u_3\}$ be the third neighbor of $u_2$,
and let $v_3$ and $v_4$ be defined similarly.
Let $G_1=G+(u_1,u_2,u_3,u_4)$, let $k_1,l_1$ be the
new vertices of $G_1$, let $k_2,l_2$ be the new vertices of 
$G_1+(u_2,u_3,u_4,u_5)$,
and let $\eta$ be the restriction
of $\eta_2$ to $G$. Let $\zeta_1$ be obtained from $\eta$
by rerouting $\eta_2(u_2u_3)$ along $\eta_2(k_1l_1)$. By
considering the path $\eta_2(l_2k_2)\cup\eta_2(k_2u_2)$ we can
extend $\zeta_1$ to a homeomorphic embedding $G+(v_2,u_2,u_4,u_5)\emb H$.
We deduce that if $G+(v_2,u_2,u_4,u_5)$ is \cfc, then the lemma
holds. Thus we may assume that that is not the case, and hence
$v_2$ and $u_5$ are adjacent in $G$ by \refthm{2.2}. 

Let $G'=G+(u_3,v_3,u_5,v_2)$ and $\eta'$ be obtained from $\zeta_1$
by first rerouting $\eta_2(v_2u_5)$ along $\eta_2(k_2l_2)\cup
\eta_2(k_2u_2)$, and then by routing the new edge along 
$\eta_2(k_2u_3)$. Since
$G'$ is \cfc\ by \refthm{2.4}, the lemma follows. \qed

\newthm{7.2}Let $G$ be a \cfc\ cubic graph, and let 
$\ldots,u_{-1},u_0,\allowbreak u_1,\ldots$ and $\ldots,v_{-1},v_0,v_1,\ldots$
be two doubly infinite sequences of (not necessarily
distinct) vertices of $G$ such that
for all integers $i$, the neighbors of $u_i$ are $u_{i-1}$, 
$u_{i+1}$ and $v_i$, and the neighbors of $v_i$ are $v_{i-2}$,
$v_{i+2}$ and $u_i$. Then there exists an integer $p\ge5$
($p\ge10$ if $p$ is even)
such that $u_i=u_{i+p}$ and $v_i=v_{i+p}$ for all integers $i$, and the vertices
$u_1,u_2,\ldots,u_p,v_1,v_2,\ldots,v_p$ are pairwise distinct.
Thus  $G$ is a biladder.

\proof 
Choose $p>0$ minimum such that for some integer $i$, one 
of $u_i,v_i$ equals one of $u_{i+p}, v_{i+p}$. Suppose first that $u_i = v_{i+p}$. 
Then $p>2$, and the neighborhood set of $u_i$ equals the neighborhood set of $v_{i+p}$, 
so one of $u_{i-1}, u_{i+1}, v_i$ equals $v_{i+p-2}$,
contrary to the choice of $p$.
If $v_i = u_{i+p}$, then similarly one of $u_{i+p-1}, u_{i+p+1}, v_{i+p}$ equals $v_{i+2}$,
again contrary to the choice of $p$.

So either $u_i = u_{i+p}$ or $v_i = v_{i+p}$; and then as before, it follows
that $u_i = u_{i+p}$ and $v_i = v_{i+p}$ for all integers $i$. It follows from
the choice of $p$ that the vertices
$u_1,u_2,\ldots,u_p,\allowbreak v_1,v_2,\ldots,v_p$ are pairwise distinct.~\qed

\newthm{7.3}Let $G,H$ be \cfc\ cubic graphs, let $H$ be a long $2$-extension
of $G$, 
and assume that there does not exist
a handle expansion $G'$ of $G$ which admits a \he\ 
$G'\emb H$. 
Then both $G,H$ are biladders.

\proof 
Since $H$ is a $2$-extension of $G$, there exist vertices $v_1,v_3,u_3,u_2$ of $G$ and 
vertices $a_1,a_2,a_3,a_4$ of $G_1=G+(v_1,v_3,u_3,u_2)$ such that 
$H = G_1+(a_1,a_2,a_3,a_4)$.
Let  $k_1,l_1$ be the new vertices of $G_1$. 
Then $G_1$ is not \cfc,
and so by \refthm{1extquadconn} we may assume that $v_3$ is adjacent to $u_3$ in $G$.
Thus $G_1$ is \qc\ and has a unique quadrangle $C_1$, where $C_1$ has
vertex-set $v_3,u_3,l_1,k_1$.
Furthermore, one of $\{a_1,a_2\}$,
$\{a_3,a_4\}$ is equal to one of $\{v_3,u_3\}$, $\{u_3,l_1\}$, $\{l_1,k_1\}$
or $\{k_1,v_3\}$. From the symmetry (and  making use of the \he\ obtained from
the canonical \he\ $G\emb H$ by rerouting $v_3u_3$ along $k_1l_1$) we may assume
that either $\{a_1,a_2\}=\{v_3,u_3\}$, or $\{a_1,a_2\}=\{u_3,l_1\}$.
Let $v_2\ne u_3$ be a neighbor of $u_2$ in $G$. In the former case,
since $G+(a_1,a_2,a_3,a_4)$ is not \cfc, we may assume from the
symmetry that $\{a_3,a_4\}=\{u_2,v_2\}$, in which case we obtain
a contradiction from \refthm{7.1} applied to the path of $G$ with vertex-set
$\{v_1,v_3,u_3,u_2,v_2\}$.

We may therefore assume that  $\{a_1,a_2\}=\{u_3,l_1\}$, and further
(by replacing $v_2$ if necessary)
that  $\{a_3,a_4\}=\{v_2,v_4\}$, where $v_4\ne u_2$ is a neighbor of $v_2$. 
Then $v_1\ne v_4$, because $G$ is \cfc.
Thus $G$ has a path $P$  with vertex-set 
$v_1,v_3,u_3,u_2,v_2,v_4$
(in order) such that $H=G\&(v_1,v_3,u_3,u_2,v_2,v_4)$.
(The $\&$ operator was defined prior to~\refthm{2.3}).
Let $u_1$ be the neighbor of $u_2$ not on $P$, and let
$u_4$ be the neighbor of $u_3$ not on $P$. Assume that for some
integers $m,n$ with $m\le1$ and $n\ge4$
we have already constructed (not necessarily distinct) vertices $u_m,u_{m+1},\ldots,u_n$, 
$v_m,v_{m+1},\ldots,v_n$ of $G$ such that for all 
$i=m+1,m+2,\ldots,n-1$
\item{(i)}$u_i$ is adjacent in $G$ to $u_{i+1}$ and $u_m$ is 
adjacent in $G$ to $u_{m+1}$,
\item{(ii)}$u_i$ is adjacent in $G$ to $v_i$,
\item{(iii)}$v_{i-1}$ is adjacent in $G$ to $v_{i+1}$,
\item{(iv)}there exists a homeomorphic embedding
$$\eta_n:G\&(v_{n-3},v_{n-1},u_{n-1},\allowbreak u_{n-2},
v_{n-2},v_n)\emb H,$$

\noindent
and
\item{(v)}there exists a homeomorphic embedding
$$\eta_m:G\&(v_m,v_{m+2},u_{m+2},\allowbreak u_{m+1},
v_{m+1},v_{m+3})\emb H.$$ 

\noindent We shall construct $u_{m-1},v_{m-1},u_{n+1},v_{n+1},
\eta_{m-1},\eta_{n+1}$
such that (i)-(v) are satisfied for all $i=m,m+1,\ldots,n$.

Let $L=G\&(v_{n-3},v_{n-1},u_{n-1},u_{n-2},v_{n-2},v_n)$, and
let $k,l,k',l'$ be the new vertices of $L$. Let $\eta'$ be
obtained from the restriction of $\eta_n$ to $G$ by rerouting
$\eta_n(v_{n-1}u_{n-1})$ along $\eta_n(kl)$. By considering the
path $\eta_n(k'l')$ we can extend $\eta'$ to a homeomorphic
embedding $\eta'':G+(u_{n-1},u_n,v_{n-2},v_n)\emb H$.
Since $G+(u_{n-1},u_n,v_{n-2},v_n)$ is not \cfc\ by
hypothesis, we deduce from \refthm{2.2} that $u_n$, $v_n$ are adjacent.
Let $u_{n+1}\not\in\{u_{n-1},v_n\}$ be the third neighbor of $u_n$,
and let $v_{n+1}\not\in\{u_{n-1},v_{n-3}\}$ be the third neighbor
of $v_{n-1}$.
By considering the \he\ $\eta''$ and the path $\eta_n(u_{n-1}v_{n-1})$ 
we can construct a \he\ 
$\eta_{n+1}:G\&(v_{n-2},v_{n},u_{n},u_{n-1},v_{n-1},
\allowbreak v_{n+1})\emb H$.
The vertices $u_{m-1},v_{m-1}$ and homeomorphic
embedding $\eta_{m-1}$ are defined analogously.

This completes the definition of two doubly infinite sequences
of vertices $\ldots u_{-1},u_0,\allowbreak u_1,\ldots$ and 
$\ldots v_{-1},v_0,v_1,\ldots$ of $G$ such that (i), (ii), (iii) hold
for all integers $i$. It follows from \refthm{7.2} that both $G,H$ are
biladders, as required. \qed

\newthm{7.4}Let $G$, $G_1$ be biladders, where
$|V(G_1)|=|V(G)|+4$ and $|V(G)|\not\in\{10,20\}$, and let 
$G_2$ be a  handle expansion of $G_1$.
Then there exist a handle expansion $G'$ of $G$
and a \he\ $G'\emb G_2$.

\proof Let us assume that the vertices of  $G_1$ are numbered
$u_0,u_1,\ldots,u_{p+1},v_0,v_1,\ldots,v_{p-1}$,
as in the definition of biladder.
The edges of the form $u_iv_i$ will be called \dfn{rungs}.
Let us say that two edges $e,f$ in a graph are \dfn{diverse} if they
share no end and no end of $e$ is adjacent to an end of $f$.
It follows by inspection that 
if $e,f$ are two diverse edges of $G_1$,
then there exist two consecutive rungs such that they are not equal to
$e,f$ and upon the deletion of the rungs
and suppression of the resulting degree two vertices
the edges (corresponding to) $e,f$  remain diverse in the smaller biladder.
%
Since deleting two consecutive rungs and suppressing the resulting
degree two vertices produces a graph isomorphic to $G$,
we deduce that the theorem holds.~\qed


The following variation of \refthm{7.4} is easy to see.

\newthm{7.5}Let $G$, $G_1$ be biladders, where $|V(G_1)|=|V(G)|+4$, and
let $G_2$ be a circuit expansion of $G_1$.
Then there exist a circuit expansion
$G'$ of $G$ and a homeomorphic embedding $\eta':G'\emb G_2$.

The following theorem implies (1.3) and (1.4).

\newthm{7.6}Let $G,H$ be non-isomorphic \cfc\ cubic graphs, 
assume that $H$ topologically contains $G$,
and assume that not both $G,H$ are biladders. 
Assume further that if 
$G$ is isomorphic
to the Petersen graph, then $H$ does not topologically contain the
biladder on $14$ vertices, and if
$G$ is isomorphic
to the Dodecahedron, then $H$ does not topologically contain the
biladder on $24$ vertices.
Then there exist a \cfc\  cubic graph 
$G'$ and a homeomorphic embedding $G'\emb H$ such that 
$G'$ is either a handle  or circuit expansion of $G$.
Moreover, if $H$ is \dc, then $G'$ can be chosen to be a handle
expansion.

\proof We proceed by induction on $|V(H)|-|V(G)|$. Let $G,H$
be as stated, and assume that the theorem holds for all pairs
$G',H'$ with $|V(H')|-|V(G')|<|V(H)|-|V(G)|$. By \refthm{6.2}
there exist a \cfc\ cubic graph $G_1$ and a homeomorphic
embedding $G_1\emb H$ such that $G_1$ is a $1$- or
$2$-extension or a circuit expansion of $G$.
If $H$ is \dc, then by \refthm{6.1} $G_1$ can be chosen
to be a $1$- or $2$-extension of $G$. We may assume that
$G_1$ is a $2$-extension of $G$, for otherwise the conclusion
of the theorem is satisfied.
 From \refthm{7.3} we deduce that either the conclusion of the theorem 
is satisfied, or  both $G,G_1$ are biladders, and so we may assume 
the latter. Thus $|V(G)|\ne20$ by the hypothesis of the theorem.
By the induction hypothesis applied
to the pair $G_1,H$ we deduce that there exist a
handle or circuit expansion $G_2$ of $G_1$ and a homeomorphic
embedding $G_2\emb H$. Moreover,
if $H$ is \dc, $G_2$ can be chosen to be a handle expansion.
By \refthm{7.4} and \refthm{7.5} there exist a handle or
circuit expansion $G'$ of $G$ and a homeomorphic embedding 
$\eta':G'\emb H$. 
Moreover, if $H$ is \dc, then $G'$ is a handle expansion. 
Thus the pair $G',\eta'$ satisfies the conclusion of the theorem. \qed

\beginsection ACKNOWLEDGEMENT

We thank Daniel P. Sanders for carefully reading the manuscript, and
for providing helpful comments.

\beginsection REFERENCES

\baselineskip 14pt

\def\JCTB{{\it J.~Combin.\ Theory Ser.\ B}}

\bibitem{AldHolJac}E. R. L. Aldred, D. A. Holton, B. Jackson, 
Uniform cyclic edge connectivity in cubic graphs,
{\it Combinatorica} {\bf 11} (1991), 81--96.

\bibitem{Bar}D. Barnette, On generating planar graphs,  
{\it Discrete Math.} {\bf 7} (1974), 199--208.

\bibitem{But}J. W. Butler, 
A generation procedure for the simple $3$--polytopes with cyclically 
$5$--connected graphs,
{\it Canad. J. Math.} {\bf 26} (1974), 686--708.

\bibitem{DieGT} R.~Diestel, Graph Theory,
Springer, 2010.

\bibitem{EdwSanSeyTho} K.~Edwards, D.~P.~Sanders, P.~D.~Seymour and R.~Thomas,
Three-edge-colouring doublecross cubic graphs, 
\JCTB\ {\bf119} (2016), 66--95.

\bibitem{Kur}K. Kuratowski, Sur le probl\`eme des courbes 
gauches en topologie, {\it Fund. Math.} {\bf 15} (1930), 271--283.

\bibitem{McCPhD}W.~McCuaig, 
Edge-reductions in cyclically $k$--connected cubic graphs,
Ph.~D. thesis, University of Waterloo, Waterloo, Ontario, October 1987.
 
\bibitem{McCEdge}W.~McCuaig, 
Edge-reductions in cyclically $k$--connected cubic graphs,
{\it J. Combin. Theory Ser. B} {\bf 56} (1992), 16--44.

\bibitem{RobSanSeyTho4CT}
N.~Robertson, D.~P.~Sanders, P.~D.~Seymour and R.~Thomas,
The four-colour theorem, {\it\JCTB} {\bf 70} (1997), 2-44.

\bibitem{RobSeyGM9}N.~Robertson and P.~D.~Seymour, 
Graph Minors IX. Disjoint crossed paths, 
\JCTB\ {\bf 49} (1990), 40--77.

\bibitem{RobSeyThoTut} N.~Robertson, P.~D.~Seymour and R.~Thomas,
Tutte's edge-coloring conjecture,
\JCTB\ {\bf 70} (1997), 166--183.

\bibitem{RobSeyThoCubic}N.~Robertson, P.~D.~Seymour and R.~Thomas,
Excluded minors in cubic graphs, {\tt arXiv:1403.2118}.

\bibitem{TutConvex}W. T. Tutte, Convex representations of graphs,  
{\it Proc. London Math. Soc.} {\bf 10} (1960), 304--320.

\bibitem{TutGeom}W.~T.~Tutte, A geometrical version of the Four Colour Problem, ``Combinatorial Mathematics and its Applications", Bose and Dowling Eds.,  
{\it The University of North Carolina Press} (1969), 553-560.

\vfill

\noindent
This material is based upon work supported by the National Science Foundation.
Any opinions, findings, and conclusions or recommendations expressed in
this material are those of the authors and do not necessarily reflect
the views of the National Science Foundation.
\eject

\end